\documentclass[sn-mathphys]{sn-jnl}

\jyear{2021}%
\usepackage{amsfonts}
\usepackage{bbm}
\usepackage{amssymb,amsfonts,boxedminipage}
\usepackage{amsthm}
\usepackage{epsfig,graphicx,picinpar,subfigure, epstopdf}  

\newtheorem{thm}{Theorem}[section]
\newtheorem{exmp}[thm]{Example}

\newcommand{\abs}[1]{\left\vert#1\right\vert}

\newcommand{\simiid}{\stackrel{\mathrm{iid}}\sim}

\newcommand{\mrd}{\mathrm{d}}
\newcommand{\mbe}{\mathbb{E}}

\newcommand{\yobs}{y_{\mathrm{obs}}}
\newcommand{\pabc}{\tilde{p}}

\newcommand{\mbp}{\mathbb{P}}

\usepackage{natbib}

\begin{document}

\title[An adaptive MPMC method for likelihood-free inference]{An adaptive mixture-population Monte Carlo method for likelihood-free inference}


\author[1]{\fnm{Zhijian} \sur{He}}\email{hezhijian@scut.edu.cn}

\author[1]{\fnm{Shifeng} \sur{Huo}}\email{18022320693@163.com}

\author*[2]{\fnm{Tianhui} \sur{Yang}}\email{yth.tianhui@gmail.com}

\affil[1]{\orgdiv{School of Mathematics}, \orgname{South China University of Technology}, \orgaddress{\city{Guangzhou}, \postcode{510641},  \country{China}}}


\affil[2]{\orgdiv{Department of Mathematics}, \orgname{Harbin Institute of Technology at Weihai}, \orgaddress{\city{Weihai}, \postcode{264209},  \country{China}}}


\abstract{
	This paper focuses on variational inference with intractable likelihood functions that can be unbiasedly estimated. A flexible variational approximation based on Gaussian mixtures is developed, by adopting the mixture population Monte Carlo (MPMC) algorithm in \cite{cappe2008adaptive}. 
	MPMC updates iteratively the parameters of mixture distributions with importance sampling computations, instead of the complicated gradient estimation of the optimization objective in usual variational Bayes. Noticing that MPMC uses a fixed number of mixture components, which is difficult to predict for real applications, we further propose an automatic component--updating procedure to derive appropriate number of components.
	The derived adaptive MPMC algorithm is capable of finding good approximations of the multi-modal posterior distributions even with a standard Gaussian as the initial distribution, as demonstrated in our numerical experiments.
\\
\\

}




\maketitle

\section{Introduction}
In this paper, we are interested in variational inference (VI) with the aim of finding a tractable distribution $q(\theta)$ that ``best" approximates the posterior distribution $\pi(\theta)=p(\theta \vert  \yobs)\propto p(\theta)p(\yobs\vert \theta)$, where $p(\theta)$ is the prior and $p(\yobs\vert \theta)$ is the likelihood function of the observed data $\yobs$.  VI is carried out by finding the member from some chosen parametric family $\mathcal{Q}$ to minimize specific divergence between probability distributions. It is usually computationally effective compared to simulation-based methods, such as Markov chain Monte Carlo (MCMC). The accuracy of VI depends on the range of the class $\mathcal{Q}$ and the divergence metric. Kullback--Leibler (KL) divergence, also known as relative entropy, is widely used in VI, where the KL divergence between probability densities $q_1$ and $q_2$ is defined as
\begin{equation*}
\mathrm{KL}(q_1 \Vert  q_2)=\int \log\left(\frac{q_1(\theta)}{q_2(\theta)}\right)q_1(\theta)\mrd\theta.
\end{equation*}
It should be noted that KL is asymmetric in replacing $q_1$ and $q_2$. Between the two forms, many variational methods (e.g., variational Bayes) minimize the form of $\mathrm{KL}(q \Vert \pi)$, whereas the expectation propagation (EP) algorithm \citep{Min:2013} considers the divergence $\mathrm{KL}(\pi \Vert  q)$. 
See \cite{Min:2005} for a good discussion of the relative merits of each divergence.

In many applications typically arising from natural and social sciences \citep{Beau:2010,davi:2012}, however, the likelihood functions $p(\yobs\vert \theta)$ concerning probabilistic models are intractable. 
In state space-space models \citep{durbin:2012}, the likelihood function is an intractable high-dimensional integral over the state variables governed by a Markov process. Approximate Bayesian computation (ABC) introduced in \cite{tavare1997} is a popular approach for Bayesian inference
when the likelihood function is analytically or computationally intractable but it is possible to simulate data from the model. Many ABC algorithms are based on sampling the ABC posteriors, including acceptance-rejection sampling, importance sampling (IS), and MCMC \citep{mar:2003}.
We refer to \cite{sisson:2018} and references therein for a comprehensive review of ABC. 

Implementing VI is difficult in the likelihood-free setting. In the context of ABC, due to the intractability of ABC posteriors, traditional VI algorithms fail to work. Some adaptations have been made to the likelihood-free setting. EP-ABC, first introduced in \cite{BC:2014}, is an adaption of EP to the ABC setting. However, EP-ABC works only for likelihoods with a special structure, and is not designed for multimodal distributions.  As shown in \cite{BC:2014}, when the posterior is multimodal, EP-ABC tends to produce an approximation covering all the modes, hence leading to a larger support. A  variational Bayes with intractable likelihood (VBIL) method is developed in \cite{tran:2017}, which can be applied to commonly used statistical models without requiring an analytical solution to model-based expectations. However, VBIL aims at minimizing the KL of densities on an extended space  $\mathrm{KL}(q(\theta,z) \Vert  \pi(\theta,z))$, where $\pi(\theta,z)$ admits the posterior $\pi(\theta)$ as its marginal and $z$ as a latent variable. This differs from the usual way of minimizing $\mathrm{KL}(q(\theta) \Vert  \pi(\theta))$. To recover VI based on $\mathrm{KL}(q(\theta) \Vert \pi(\theta))$ in the likelihood-free setting, \cite{he:2021} uses multilevel Monte Carlo \citep{Giles2015,Rhee2015} to derive unbiased estimates of the gradient of the KL. 
On the other hand, \cite{ong:2018} modified the VBIL method to work with unbiased log-likelihood estimates in the synthetic likelihood framework, resulting in the variational Bayes  synthetic likelihood (VBSL) method. A limitation of VBSL is its assumption on the normality of summary statistics.
Many variational methods including VBIL and VBSL are based on stochastic gradient decent (SGD) algorithms, which require an unbiased estimator  of the gradient of KL.  A difficulty with SGD is that plain Monte Carlo  sampling to estimate the gradient can be error prone or inefficient. Some variance reduction methods have been adopted to improve SGD \citep[see][]{PBJ:2012,MF:2017}. 

The variational methods in the likelihood-free setting are generally not designed for posteriors with  multimodal structures. Thus they may not render a good approximation of complex posteriors. Considering this, we introduce an adaptation of the mixture population Monte Carlo (MPMC) method proposed by \cite{cappe2008adaptive} to the likelihood-free setting. MPMC finds the optimal $q(\theta)$ by minimizing $\mathrm{KL}(\pi \Vert  q)$
within a class of mixture distributions, such as Gaussian mixtures. The proposed MPMC algorithm in this paper is able to fit posteriors with multimodal and heavy-tailed structures by using a proper mixture distribution. Additionally, MPMC updates iteratively the parameters of  mixture distributions by IS computations rather than SGD algorithms. In other words, MPMC does not have to calculate the gradient of KL, which is difficult in the likelihood-free setting. In the original form of the MPMC method \citep{cappe2008adaptive}, the number of mixture components $D$ is fixed throughout. For practical problems, how to choose a proper $D$ in advance is still not clear.
Many mixture parameters need to be initialized when $D$ is large.
To handle these issues, we propose an adaptive procedure to update $D$ automatically. Meanwhile, we provide a strategy to initialize a new mixture component when updating the  number of mixture components. The whole process of our proposed method is automatic with minimal inputs.
The main contributions of this paper are:
\begin{itemize}
	\item We propose a flexible VI method to approximate posteriors with intractable likelihood function, following a way instead of the SGD paradigm.
	\item An adaptive algorithm is designed to guide the components in the mixture distribution, where a new component would be added for a stronger approximating ability while the unnecessary ones would be ruled out from the mixture proposal.
	\item Numerical experiments demonstrate the effectiveness and efficiency of the proposed method.
\end{itemize}

The rest of this paper is organized as follows.  Section~\ref{sec:mpmc} formulates the problems of likelihood-free inference. In Section~\ref{sec:main}, we first review the MPMC method and then  adapt it to the likelihood-free setting. We further show how to initialize a new mixture component and propose some stopping rules. In Section~\ref{sec:num}, we provide some numerical examples. Section~\ref{sec:con} concludes this paper.

\section{Problem formulation}\label{sec:mpmc}
Let $\theta\in\Theta\subset\mathbb{R}^p$ be the parameter of interest with prior density $p(\theta)$, $\yobs$ be the observed data and $p(\yobs\vert \theta)$ be the likelihood. Our target is to approximate the posterior distribution 
\begin{equation}
\pi(\theta):=p(\theta\vert \yobs)= \frac{p(\theta)p(\yobs\vert \theta)}{p(\yobs)},\label{eq:post}
\end{equation}
where $p(\yobs) = \int p(\theta)p(\yobs\vert \theta) \mrd \theta$ is usually an  unknown constant (called the marginal likelihood or evidence).
In many applications, however, the likelihood $p(\yobs\vert \theta)$ is analytically intractable, such as the two examples given below. 

\begin{exmp}[Approximate Bayesian computation]\label{exam:ABC}
	ABC is a generic method in likelihood-free inference provided that it is easy to generate $y\sim  p (y \vert  \theta)$. However, ABC methods do not target the exact posterior, but an approximation to some extent. More specially, let $S(\cdot):\mathbb{R}^n\to \mathbb{R}^d$ be a vector of summary statistics, and $K_h(\cdot,\cdot)$ be  a $d$-dimensional kernel density with bandwidth $h>0$. 
	ABC posterior density of  $\theta$ is given by
	\begin{equation}
	p_{\mathrm{ABC}}(\theta\vert \yobs)\propto p(\theta)\pabc(\yobs\vert \theta),
	\end{equation}
	where the intractable likelihood is given by
	\begin{equation}
	\begin{aligned}
	\pabc(\yobs\vert \theta)&=\int K_h(S(\yobs),S(y))p(y\vert \theta) \mrd y\\&=\mathbb{E}_{p(y\vert \theta)}[K_h(S(\yobs),S(y))].
	\end{aligned}\label{eq:abc}
	\end{equation}
	If $S(y)$ is a sufficient statistic, then $p_{\mathrm{ABC}}(\theta\vert \yobs)$ converges to the exact posterior $p (\theta \vert \yobs)$ as $h\to 0$. Otherwise, $p_{\mathrm{ABC}}(\theta\vert \yobs)$ converges to the posterior $p (\theta \vert  S(\yobs))$ as $h\to 0$. There is a gap between $p (\theta \vert  S(y))$ and $p (\theta \vert  y)$. We refer to \cite{blum_comparative_2013} for a comparative review of choosing low-dimensional summary statistics.
\end{exmp}

\begin{exmp}[State-space models]
	In state-space models, the likelihood is a high-dimensional integral over the state variables governed by a Markov process, which is usually either computationally prohibitive,
	or simply not possible. Suppose that the observations $y_t$ are observed in time order. At time $t$, the distribution of $y_t$ conditional on a state (latent) variable $x_t$ is independently distributed as $y_t \vert  x_t\sim g_t(y_t \vert  x_t,\theta)$, and the state variables $\{x_t\}_{t\ge 1}$ are Markov chain with $x_1\sim \mu_{\theta}(x_1)$ and $x_t \vert  x_{t-1}\sim f(x_t\vert   x_{t-1},\theta)$. The likelihood of the data $\yobs = (y_1,\dots,y_T)$ is 
	\begin{equation}
	p(\yobs\vert \theta)=\int p(\yobs\vert  x,\theta)p(x\vert \theta)\mrd x = \mbe_{p(x \vert  \theta)}[p(\yobs \vert  x,\theta)],\label{eq:ssm}
	\end{equation}
	where $x= (x_1,\dots,x_T)$, 
	$$p(x\vert \theta) = \mu_{\theta}(x_1)\prod_{t=2}^T f_t(x_t \vert  x_{t-1},\theta) \text{~and~} p(\yobs \vert  x,\theta)= \prod_{t=1}^T g_t(y_t, \vert  x_t,\theta).$$
\end{exmp}

In this paper, we consider the setting in which the likelihood $p(\yobs\vert \theta)$ is intractable, but can be estimated unbiasedly. Let $\Psi(X;\yobs,\theta)$ be an unbiased estimator of  $p(\yobs\vert \theta)$, where $X\sim p(x \vert  \theta)$ encapsulates all the random elements in estimating the likelihood.  We thus arrive at a model-based expectation  
\begin{equation}\label{eq:setting}
p(\yobs\vert \theta) = \mbe_{X}[\Psi(X;\yobs,\theta) \vert  \theta],
\end{equation}
which does not have an analytical solution. It is clear that the two examples above fulfill the form \eqref{eq:setting}. More concrete examples can be found in Section~\ref{sec:num}. In the following, we use a proper mixture distribution to fit the posterior $\pi(\theta)$ given by \eqref{eq:post}.

\section{An adaptive MPMC algorithm}\label{sec:main}
We first review the idea of the MPMC algorithm proposed in \cite{cappe2008adaptive}. MPMC takes a mixture density $$ q_{(\alpha,\eta)}(\theta)=\sum_{d=1}^D{\alpha_dq_d(\theta;\eta_d)} $$ with chosen parameters $\alpha=(\alpha_1,\dots,\alpha_D)$ and $\eta=(\eta_1,\dots,\eta_D)$ to minimize the KL divergence between the target density $\pi$ and the mixture $q_{(\alpha,\eta)}$:
\begin{equation}
\mathrm{KL}(\pi \Vert  q_{(\alpha,\eta)})=\int_{\Theta} \log\left(\frac{\pi(\theta)}{q_{(\alpha,\eta)}(\theta)}\right)\pi(\theta)\mrd\theta,\label{KL}
\end{equation}
where $\alpha$ satisfies the partition of unity constraints $\sum_{d=1}^D \alpha_d=1$ with entries $\alpha_d\ge 0$ throughout this paper.
Minimizing \eqref{KL} is equivalent to maximizing the following objective function:
\begin{equation}
L(\alpha,\eta):=\mbe_{\pi}[\log(q_{(\alpha,\eta)}(\theta))]=\int_{\Theta} \log(q_{(\alpha,\eta)}(\theta))\pi(\theta)\mrd\theta.
\end{equation}
To solve this optimization problem,  \cite{cappe2008adaptive} used the integrated expectation-maximization (EM) algorithm, which is quite similar to the original EM algorithm proposed by \cite{em:1977} with the E-step replaced by IS computations.

Let $Z$ be the component indicator (a discrete random variable) taking values in $\mathcal{D}:=\{1,\dots,D\}$ 
such that the joint density of $Z$ and $\theta$ satisfies
$$f(z,\theta) = \alpha_z q_z(\theta;\eta_z).$$
Note that $\mbp(Z=z)=\alpha_z$ for $z\in \mathcal{D}$, and
$$\int_{\mathcal{D}} f(z,\theta) \mrd z = \sum_{d=1}^D{\alpha_dq_d(\theta;\eta_d)} =q_{(\alpha,\eta)}(\theta),$$
where the integration  is carried out under the counting measure.
Let $g(z)$ be a probability mass function defined on the domain $\mathcal{D}$. Assume $g(z)>0$ for all $z\in \mathcal{D}$. By a change of measure, we have
\begin{align}
L(\alpha,\eta) &=\mbe_{\theta\sim\pi}\left[\log\left(\int_{\mathcal{D}} f(z,\theta) \mrd z\right)\right] \notag\\
&=\mbe_{\theta\sim\pi}\left[\log \mbe_{Z\sim g}[f(Z,\theta)/g(Z) \vert  \theta]\right]\notag\\
&\ge \mbe_{\theta\sim\pi}\left[ \mbe_{Z\sim g}[\log(f(Z,\theta)/g(Z)) \vert  \theta]\right]:=\ell(g,\alpha,\eta),\label{eq:ell}
\end{align}
where the last inequality is due to Jensen's inequality. Let
\begin{equation}
f(z \vert  \theta,\alpha,\eta) := \frac{f(z,\theta)}{\int_{\mathcal{D}} f(z,\theta)\mrd z}=\frac{\alpha_z q_z(\theta;\eta_z)}{\sum_{d=1}^D\alpha_d q_d(\theta;\eta_d)}
\end{equation}
be the mixture posterior probabilities.
It is clear that 
\begin{equation}
L(\alpha,\eta)=\ell(f(z \vert  \theta,\alpha,\eta),\alpha,\eta).
\end{equation}

As in the usual EM algorithm, the proposed updating mechanism is given by
\begin{align}
(\alpha^{t+1},\eta^{t+1}) &= \arg\max_{(\alpha,\eta)} \ell(f(z \vert  \theta,\alpha^t,\eta^t),\alpha,\eta)\label{eq:emiter}\\
&= \arg\max_{(\alpha,\eta)} \mbe_{\theta\sim\pi}\left[ \mbe_{Z\sim f(z \vert  \theta,\alpha^t,\eta^t)}[\log(f(Z,\theta)/f(Z \vert  \theta,\alpha^t,\eta^t)) \vert  \theta]\right]\notag\\
&=\arg\max_{(\alpha,\eta)} \mbe_{\theta\sim\pi}\left[ \mbe_{Z\sim f(z \vert  \theta,\alpha^t,\eta^t)}[\log(f(Z,\theta)) \vert  \theta]\right]\notag\\
&=\arg\max_{(\alpha,\eta)} \mbe_{\theta\sim\pi}\left[ \mbe_{Z\sim f(z \vert  \theta,\alpha^t,\eta^t)}[\log(\alpha_Zq_Z(\theta;\eta_Z)) \vert  \theta]\right].\label{eq:mstep}
\end{align}
It is easy to see that
\begin{equation}
\begin{aligned}
L(\alpha^{t+1},\eta^{t+1})&\ge \ell(f(z \vert  \theta,\alpha^t,\eta^t),\alpha^{t+1},\eta^{t+1})\\
&\ge \ell(f(z \vert  \theta,\alpha^t,\eta^t),\alpha^{t},\eta^{t})=L(\alpha^{t},\eta^{t}).
\end{aligned} \label{eq:lw1}
\end{equation}
This first inequality in \eqref{eq:lw1} is due to \eqref{eq:ell}, and the second inequality in \eqref{eq:lw1} results from \eqref{eq:emiter}. This implies that $(L(\alpha^{t},\eta^{t}))_{t\ge 1}$ is a non-decreasing sequence.

It should be noticed that the maximization step \eqref{eq:mstep} reduces to solve
\begin{align*}\alpha^{t+1}&=\arg\max_{\alpha} \mbe_{\theta\sim\pi}\left[ \mbe_{Z\sim f(z \vert  \theta,\alpha^t,\eta^t))}[\log(\alpha_Z) \vert  \theta]\right]\\&=\arg\max_{\alpha}\mbe_{\theta\sim\pi}\left[ \sum_{d=1}^D f(d \vert  \theta,\alpha^t,\eta^t)\log(\alpha_d) \right]\\
&=\arg\max_{\alpha}\sum_{d=1}^D \mbe_{\theta\sim\pi}[f(d \vert  \theta,\alpha^t,\eta^t)]\log(\alpha_d),
\end{align*}
and
\begin{align*}\eta^{t+1}&=\arg\max_{\eta} \mbe_{\theta\sim\pi}\left[ \mbe_{Z\sim f(z \vert  \theta,\alpha^t,\eta^t))}[\log(q_Z(\theta;\eta_Z)) \vert  \theta]\right]\\&=\arg\max_{\eta}\mbe_{\theta\sim\pi}\left[ \sum_{d=1}^D f(d \vert  \theta,\alpha^t,\eta^t)\log(q_d(\theta;\eta_d)) \right]\\
&=\arg\max_{\eta} \sum_{d=1}^D \mbe_{\theta\sim\pi}[f(d \vert  \theta,\alpha^t,\eta^t)\log(q_d(\theta;\eta_d)) ].
\end{align*}
As a result, for $d=1,\dots,D$, we get
\begin{align}
\alpha_d^{t+1} &= \mbe_{\theta\sim\pi}[f(d \vert  \theta,\alpha^t,\eta^t)],\label{eq:alphaIS}\\
\eta_d^{t+1} &= \arg\max_{\eta_d} \mbe_{\theta\sim\pi}[f(d \vert  \theta,\alpha^t,\eta^t)\log(q_d(\theta;\eta_d))]\label{eq:etaIS}.
\end{align}

Note that the updating mechanism involves unknown expectations with respect to the target density $\pi$. In the MPMC algorithm, the parameters of the mixture distribution at the current iteration are estimated by using the self-normalized IS with the mixture distribution at the last iteration as the proposal. Let $\hat\alpha^t = (\hat\alpha_1^t,\dots,\hat\alpha^t_D)$ and $\hat\eta^t = (\hat\eta_1^t,\dots,\hat\eta^t_D)$ be the estimated parameters at the $t$-th iteration. Consider a general problem of estimating an expectation $I(\varphi):=\mbe_\pi[\varphi(\theta)]$. Let $ (\theta_{i,t})_{1\leq i\leq N} $ be the iid sample generated from the mixture $ q_{(\hat\alpha^t,\hat\eta^t)}$ at the $t$-th iteration, and let $\pi_u$ be an unnormalized version of $\pi$. In the context of Bayesian inference, we take $\pi_u(\theta)=p(\theta)p(\yobs\vert \theta)$.   Denote $w(\theta)=\pi_u(\theta)/q_{(\hat\alpha^t,\hat\eta^t)}(\theta)$. The self-normalized IS is carried out by noticing  $$I(\varphi)=\frac{\mbe_{q_{(\hat\alpha^t,\hat\eta^t)}}[w(\theta)\varphi(\theta)]}{\mbe_{q_{(\hat\alpha^t,\hat\eta^t)}}[w(\theta)]},$$ yielding a ratio estimator
$$ \hat I_N(\varphi)=\frac{\sum_{i=1}^{N}w(\theta_{i,t})\varphi(\theta_{i,t})}{\sum_{i=1}^{N}w(\theta_{i,t})}=\sum_{i=1}^{N}w_{i,t}\varphi(\theta_{i,t}),$$
where the normalized importance weights are
$$w_{i,t} = \frac{w(\theta_{i,t})}{\sum_{n=1}^N w(\theta_{n,t})}=\frac{\pi_u(\theta_{i,t})}{\sum_{d=1}^D \hat\alpha_d^t q_d(\theta_{i,t};\hat\eta_d^t)}\bigg/\sum_{n=1}^{N}\frac{\pi_u(\theta_{n,t})}{\sum_{d=1}^D \hat\alpha_d^t q_d(\theta_{n,t};\hat\eta_d^t)}.$$
The updating mechanisms for \eqref{eq:alphaIS} and \eqref{eq:etaIS} are then replaced by
\begin{align}
\hat\alpha_d^{t+1} &= \hat I_N(f(d \vert  \cdot,\hat\alpha^t,\hat\eta^t))=\sum_{i=1}^N w_{i,t}f(d \vert  \theta_{i,t},\hat\alpha^t,\hat\eta^t)\text{ and }\label{eq:alphaIS2}\\
\hat\eta_d^{t+1} &= \arg\max_{\eta_d} \hat I_N(f(d \vert  \cdot,\hat\alpha^t,\hat\eta^t)\log(q_d(\cdot;\hat\eta_d)))\notag\\
&=\arg\max_{\eta_d} \sum_{i=1}^N w_{i,t}f(d \vert  \theta_{i,t},\hat\alpha^t,\hat\eta^t)\log(q_d(\theta_{i,t};\hat\eta_d))\label{eq:etaIS2},
\end{align}
respectively. Also, the objective function $L(\alpha,\eta)$ at the $t$-th iteration can be estimated by
\begin{equation}
\hat{L}(\hat\alpha^{t},\hat\eta^{t}) = \hat I_N(\log (q_{(\hat\alpha^{t},\hat\eta^{t})}(\cdot))) = \sum_{i=1}^N w_{i,t} \log \left( \sum_{d=1}^D \hat\alpha_d^t q_d(\theta_{i,t};\hat\eta_d^t)\right).\label{eq:objective}
\end{equation} 

In this paper, we focus on Gaussian mixture families in which $q_d(\theta;\eta_d)$ is the density of $p$-dimensional Gaussian $N(\mu_d,\Sigma_d)$ with $\eta_d = (\mu_d,\Sigma_d)$. For this case, \cite{cappe2008adaptive} showed that the infimum for  \eqref{eq:etaIS} is reached when
\begin{align*}
\mu_d^{t+1}&=\mbe_{\theta\sim\pi}[f(d \vert  \theta,\alpha^t,\eta^t)\theta]/\alpha_d^{t+1},\\
\Sigma_d^{t+1}&=\mbe_{\theta\sim\pi}[f(d \vert  \theta,\alpha^t,\eta^t)(\theta-\mu_d^{t+1})(\theta-\mu_d^{t+1})^\top]/\alpha_d^{t+1}.
\end{align*}
Under the self-normalized IS computations, the iteration for $\eta_d=(\mu_d,\Sigma_d)$ turns out to be
\begin{align}
\hat\mu_d^{t+1}&=\frac{\sum_{i=1}^N w_{i,t}f(d \vert  \theta_{i,t},\hat\alpha^t,\hat\eta^t)\theta_{i,t}}{ \hat\alpha_d^{t+1}},\\
\hat\Sigma_d^{t+1}&=\frac{\sum_{i=1}^N w_{i,t}f(d \vert  \theta_{i,t},\hat\alpha^t,\hat\eta^t)(\theta_{i,t}-\hat\mu_d^{t+1})(\theta_{i,t}-\hat\mu_d^{t+1})^\top}{\hat \alpha_d^{t+1}},
\end{align}
where $\hat \alpha_d^{t+1}$ is given by \eqref{eq:alphaIS2}.

\subsection{Adapting MPMC for likelihood-free inference}
Now we are ready to make an adaptation of the MPMC algorithm to the intractable likelihood \eqref{eq:setting} (named as MPMC-IL). The target density is then
$$\pi(\theta)=p(\theta \vert \yobs)=\frac{p(\theta)\mbe_X[\Psi(X;\yobs,\theta) \vert  \theta]}{p(\yobs)}.$$
To facilitate analysis, we assume that the likelihood estimator $\Psi(X;\yobs,\theta)\ge 0$, which holds for our presented examples. Otherwise, we may use the decomposition $\Psi(X;\yobs,\theta) = \Psi^+(X;\yobs,\theta)-\Psi^-(X;\yobs,\theta)$ and then work on the two parts $\Psi^{\pm}$ separately, where $\Psi^+ = \max(\Psi,0)\ge 0$ and $\Psi^- = \max(-\Psi,0)\ge 0$. Define the following density on the  extended space $\Theta\times \mathcal{X}$
$$\pi_{\mathrm{ext}}(\theta,x)=\frac{p(\theta)\Psi(x;\yobs,\theta)p(x \vert  \theta)}{p(\yobs)},$$
which admits the posterior of interest $\pi(\theta)$ as its marginal.
The key in  MPMC is to compute posterior expectations of the form $$I(\varphi)=\mbe_{\pi}[\varphi(\theta)]=\mbe_{\pi_{\mathrm{ext}}}[\varphi(\theta)],$$
whose value can be estimated by the self-normalized IS. Consider the IS density on the  extended space $\Theta\times \mathcal{X}$
$$\tilde{\pi}_{\mathrm{ext}}(\theta,x)= q_{(\hat\alpha^t,\hat\eta^t)}(\theta)p(x \vert  \theta).$$
We thus rewrite the posterior expectation $I(\varphi)$ as a ratio of two expectations $$I(\varphi)=\frac{\mbe_{\tilde{\pi}_{\mathrm{ext}}}[w(\theta,X)\varphi(\theta)]}{\mbe_{\tilde{\pi}_{\mathrm{ext}}}[w(\theta,X)]},$$
where the (unnormalized) likelihood ratio 
$$ w(\theta,x):=\frac{p(\theta)\Psi(x;\yobs,\theta)}{q_{(\hat\alpha^t,\hat\eta^t)}(\theta)}\propto \frac{\pi_{\mathrm{ext}}(\theta,x)}{\tilde{\pi}_{\mathrm{ext}}(\theta,x)}.$$ 
We should note that it is easy to generate samples from the IS density $\tilde{\pi}_{\mathrm{ext}}(\theta,x)$ by using $\theta \sim q_{(\hat\alpha^t,\hat\eta^t)}(\theta)$ and $x\sim p(x \vert \theta)$.  Let $(\theta_{i,t},X_{i,t})\simiid\tilde{\pi}_{\mathrm{ext}}(\theta,x)$ for $i=1,\dots,N$ be samples generated at the $t$-th iteration. Applying the normalized IS method gives a ratio estimator
\begin{equation}
\hat I_N(\varphi)=\frac{\sum_{i=1}^{N}w(\theta_{i,t},X_{i,t})\varphi(\theta_{i,t})}{\sum_{i=1}^{N}w(\theta_{i,t},X_{i,t})}=\sum_{i=1}^{N}w_{i,t}\varphi(\theta_{i,t}),\label{eq:ratio}
\end{equation}
where 
the normalized  weights are
\begin{align}
w_{i,t} &= \frac{w(\theta_{i,t},X_{i,t})}{\sum_{n=1}^N w(\theta_{n,t},X_{i,t})}\notag\\&=\frac{p(\theta_{i,t})\Psi(X_{i,t};\yobs,\theta_{i,t})}{\sum_{d=1}^D \hat\alpha_d^t q_d(\theta_{i,t};\hat\eta_d^t)}\bigg/\sum_{n=1}^{N}\frac{p(\theta_{n,t})\Psi(X_{n,t};\yobs,\theta_{i,t})}{\sum_{d=1}^D \hat\alpha_d^t q_d(\theta_{n,t};\hat\eta_d^t)}.\label{eq:newwei}
\end{align}
The updating mechanisms are therefore obtained by using \eqref{eq:alphaIS2} and \eqref{eq:etaIS2} with the new weights \eqref{eq:newwei}. The work in
\cite{cappe2008adaptive} showed that convergence of MPMC as $N$ increases can be established using the same approach as in \cite{douc:2007a,douc:2007b}. The argument also applies for the augmented densities on $\Theta\times \mathcal{X}$ in this likelihood-free setting.

The procedure for the Gaussian mixture family with a fixed number of components $D$ is summarized in Algorithm~\ref{alg:1}, where some stopping rules will be discussed in the next subsection. To apply MPMC-IL for ABC, it suffices to set
$$\Psi(X;\yobs,\theta)=K_h(S(\yobs),S(X)),$$
in which $X$ is distributed from the likelihood $p(y \vert \theta)$.
The next subsection discusses how to optimize the number of mixture components $D$ and how to initialize a new mixture component.

\begin{algorithm}[H]
	\caption{MPMC-IL with a fixed number of mixture components}\label{alg:1}
		Initialize $\hat\alpha^0, \hat\eta^0$, $t=0$, $D$ the number of mixture components,
$N$ the size of samples at each iteration
		\begin{algorithmic}[1]
			
		
		
		\Repeat
		
		\For{$ i=1 $ to $ N $} \label{alg:1:mpmcil:2} 
		
		\State Generate $\theta_{i,t}\sim q_{(\hat\alpha^t,\hat\eta^t)}(\theta)$, $x_{i,t}\sim p(x\vert\theta_{i,t})$ \label{alg:1:ytransform}
		
		\State Set $ w_{i,t}=p(\theta_{i,t})\Psi(x_{i,t};\yobs,\theta_{i,t})/q_{(\hat\alpha^t,\hat\eta^t)}(\theta_{i,t})$
		
		\State Set $ \rho_{i,t,d}=\hat\alpha_d^tq_d(\theta_{i,t};\hat\eta_d^t)/q_{(\hat\alpha^t,\hat\eta^t)}(\theta_{i,t})$ for $d=1,\dots,D$
		\EndFor	
		
		\State Normalize $ w_{i,t}= w_{i,t}/\sum_{j=1}^{N}w_{j,t}$ for  $i=1,\dots,N$  \label{alg:1:logw}
		
		\State Estimate the parameters $\alpha$ and $\eta=(\mu,\Sigma)$ as following, for $d=1,\dots,D$, 
		\begin{align*}
		\hat\alpha_d^{t+1}&=\sum_{i=1}^{N}w_{i,t}\rho_{i,t,d}\\
		\hat\mu_d^{t+1}&=\sum_{i=1}^{N}w_{i,t}\rho_{i,t,d}\theta_{i,t}/\hat\alpha_d^{t+1}\\
		\hat\Sigma_d^{t+1}&=\sum_{i=1}^{N}w_{i,t}\rho_{i,t,d}(\theta_{i,t}- \hat\mu_d^{t+1})(\theta_{i,t}- \hat\mu_d^{t+1})^\top/\hat\alpha_d^{t+1} 
		\end{align*} \label{alg:1:mpmcil:8} 
		
		\State Compute the object function value $L_t=\hat{L}(\hat\alpha^{t},\hat\eta^{t})$ given by \eqref{eq:objective}
		
		\State $t=t+1$ 
		
		\Until some stopping rule is satisfied
	\end{algorithmic}
\end{algorithm}


\subsection{Initializing a new mixture component and stopping rules}\label{sec:adding}
To approximate the target posterior, we further assume that the number $D$ in the Gaussian mixture is not fixed and allows for adaptive changes according to the current estimation and the observations. A  Gaussian mixture with more components is expected to achieve a stronger approximate ability, thus is preferred as a good importance density. However, using a large number of components at the initial step might lead to the problem of over-fitting and make the iterations costly. Following the idea in \cite{gunawan2021flexible} which established flexible VB based a copula of a mixture of normals but for tractable likelihoods, we add components successively in the process of iterations. 

In what follows, we start with the Gaussian mixture with one component ($D=1$), then update (add or delete) its components until the MPMC-IL algorithm converges in some sense. To this end, some rules shall be specified to terminate the current MPMC-IL iterations before updating $D$.  The first and easiest way is to introduce a fixed window size $t_w$. The number of components $D$ is updated after every $t_w$ iterations. 
Instead of setting a fixed $t_w$,  on the other hand, one can terminate adaptively the current MPMC-IL iterations. 
To be more specific, we check the smoothed objective function value $\tilde L_t$, which is an average of the objective function values over a small window size to avoid the sampling noise. If the distance between adjacent smoothed function values is too small, the current proposal can be then regard as the local minimum (or approaching its limit), and the component-updating procedure is thus conducted. 
The two rules of updating $D$  are summed up below.

\begin{enumerate}[a)]
	\item Update with a fixed window size $t_w$:
	
	\begin{equation}  \label{eq:fixedwindow}
	t > t_w
	\end{equation}
	\item Update with an adaptive window size:
	\begin{equation} \label{eq:adaptwindow}
 \lvert \tilde L_t - \tilde L_{t-1} \rvert <\epsilon_{0}, 
	\end{equation} 
	where $\tilde L_t=(1/s)\sum_{l=0}^{s-1}L_{t-l}$ for $t\ge s$, otherwise $\tilde L_t=L_t$,  $L_t:=\hat{L}(\hat\alpha^{t},\hat\eta^{t})$ is the estimated objective function value at the $t$-th iteration given by \eqref{eq:objective}, and $s$ and $\epsilon_0$ are specified by the user.
\end{enumerate}

To add an effective component means to ``make up" information that has not been well exploited by the current approximation and increase the number of the mixture components heuristically. 
Suppose the current Gaussian mixture reads
\begin{equation} \label{eq:pmc:currentmodes}
q_{(\alpha,\eta)}(\theta)=\sum_{d=1}^{D} \alpha_{d} N\left(\theta \vert  {\mu}_{d}, {\Sigma}_{d}\right),
\end{equation}
where ${\alpha_{d},~{\mu}_{d}}$, and ${{\Sigma}_{d}}$ denote respectively the weight, mean, and covariance matrix of the $d$-th component delivered by the  MPMC-IL algorithm for a fixed $D$.
In order to run the MPMC-IL algorithm for the case of $D+1$ components, it is critical to  set up the initial Gaussian mixture. To this end, we take
\begin{equation} \label{eq:pmc:init}
\tilde q(\theta)=(1-\alpha_{\mathrm{add}})\times q_{(\alpha,\eta)}(\theta)+\alpha_{\mathrm{add}}\times N\left(\theta \vert  {\mu}_{\mathrm{add}}, {\Sigma}_{\mathrm{add}}\right)
\end{equation}
as its initial mixture distribution, which is a weighted average of the current mixture $q_{(\alpha,\eta)}(\theta)$ and a new mixture component $N\left(\theta \vert  {\mu}_{\mathrm{add}}, {\Sigma}_{\mathrm{add}}\right)$. In fact, the distribution \eqref{eq:pmc:init} is a Gaussian mixture with $D+1$ components with mixing weights $((1-\alpha_{\mathrm{add}})\times\alpha, \alpha_{\mathrm{add}})$. The key is how to initialize the parameters $\alpha_{\mathrm{add}}$, $\mu_{\mathrm{add}}$ and ${\Sigma}_{\mathrm{add}}$. Among which, the determination of $\mu_{\mathrm{add}}$ needs more caution, since it shall capture the information that has not been well exploited by the current approximation $q_{(\alpha,\eta)}(\theta)$. To this end,
we draw a set of samples from $\tilde{\pi}_{\mathrm{ext}}(\theta,x)=q_{(\alpha,\eta)}(\theta)p(x\vert\theta)$, and calculate their (unnormalized) likelihood ratios given by 
\begin{equation}
\mathrm{LR}(\theta,x):=\frac{p(\theta)\Psi(x;\yobs,\theta)}{q_{(\alpha,\eta)}(\theta)}\propto \frac{\pi_{\mathrm{ext}}(\theta,x)}{\tilde{\pi}_{\mathrm{ext}}(\theta,x)}.\label{eq:lr}
\end{equation}
A sample with a large likelihood ratio implies that the target distribution may not be well approximated by the current mixture around the sample point.
We therefore take the sample of $\theta$ with the largest likelihood ratio as the mean $\mu_{\mathrm{add}}$ of the new component for compensating the misspecification of the current approximation.
Moreover, in order not to destroy the current approximation achieved, the weight for the new component $\alpha_{\mathrm{add}}$  shall not be too large. The covariance matrix ${\Sigma}_{\mathrm{add}}$ can be set as  the initial value of the first component's covariance matrix or with a small scale.

%

%

Practically, we check and delete components with too small weights to ensure numerical stability before adding a new component. 
As to the global stopping rule, three factors might be considered to stop the updating process: the budget of total iterations $T_{\max}$, the maximal number of mixture components $D_{\max}$, and the threshold $\epsilon_{\mathrm{tot}}$ for the change of the (smoothed) objective function value. Hence, the global stopping rule for the whole procedure can be set when the limit $T_{\max}$ or $D_{\max}$ is achieved, or the threshold $\epsilon_{\mathrm{tot}}$ is approached. Our proposed adaptive MPMC-IL algorithm is summarized in Algorithm~\ref{algo:component--updating}. 

\begin{algorithm}
	\caption{Adaptive MPMC-IL algorithm}
	Initialize $\hat\alpha^0=1, \hat\eta^0=(\mu_1^0,\Sigma_1^0)$, $\Sigma_{\mathrm{add}}$, $N_{\mathrm{add}}$, $D=1$, $T=0$, $L_{\mathrm{last}}=0$, $\alpha_{\min}$, the sample size $N$ at each iteration, parameters for stopping rules: $t_w$, $s$, $\epsilon_0$, $T_{\max}$, $D_{\max}$, and $\epsilon_{\mathrm{tot}}$
	\begin{algorithmic}[1] 		


	\State {Run Algorithm~\ref{alg:1} under the stopping rule \eqref{eq:fixedwindow} or \eqref{eq:adaptwindow} to get outputs of the mixture distribution with parameters $\alpha$ and $\eta$, the number of iterations $\Delta T$, and the final (smoothed) objective function value $L_{\mathrm{cur}}$}\label{step1}
	
	\State Update the total number of iterations: {$T = T+\Delta T$}
	\State Set $\Delta L=\abs{L_{\mathrm{cur}}-L_{\mathrm{last}}}$ and $L_{\mathrm{last}}=L_{\mathrm{cur}}$
	\If {$\min\limits_{i\in\{1,\cdots,d\}} \alpha_i<\alpha_{\min}$} \label{alg:2:remove cond}
	\State{Set $r = \mathop{\arg\min}\limits_{i\in\{1,\cdots,d\}} \alpha_i$}
	\State{Update parameters $\alpha=\alpha/(1-\alpha_r)\backslash \{\alpha_r\}$, $\eta=\eta \backslash \{\eta_r\}$, and $D=D-1$}
	\EndIf	
	
	\State {Generate $\theta_i\sim q_{(\alpha,\eta)}(\theta)$ and $x_{i}\sim p(x\vert\theta_{i})$, $i = 1,\dots,N_\mathrm{add}$}  \label{alg:2:saddsize}
	
	\State { Set the added mean $\mu_{\mathrm{add}}=\theta_{i^*}$, where $i^*=\mathop{\arg\max}\limits_{i\in\{1,\cdots,N_{\mathrm{add}}\}} \mathrm{LR}(\theta_i,x_i)$ and $\mathrm{LR}(\theta_i,x_i)$ is given by \eqref{eq:lr}}
	
	\State { Update the initial values for the next iteration: $D=D+1$, $\hat\alpha^0 = ((1-\alpha_{\mathrm{add}})\times\alpha, \alpha_{\mathrm{add}})$, and $\hat\eta^0=(\eta,\eta_{\mathrm{add}})$ with $\eta_{\mathrm{add}}=(\mu_{\mathrm{add}},\Sigma_{\mathrm{add}})$} \label{alg:2:alpha}
	\State { Go to Step~\ref{step1}} with the updated initial values $\hat\alpha^0, \hat\eta^0$ and $D$ until the global stopping rule is satisfied.
\end{algorithmic} \label{algo:component--updating}
\end{algorithm}

\section{Numerical experiments}\label{sec:num}
In this section, we present two examples to demonstrate the effectiveness of the proposed adaptive MPMC-IL algorithm. Different aspects are examined for the g-and-k model first, then follows a test with real data to show the robustness of the proposed algorithm.

\subsection{The g-and-k distribution}
The univariate g-and-k distribution is a flexible unimodal distribution exhibiting skewness and kurtosis \cite{RM:2002}. 
Its density function has no closed form, but its quantile function takes the form of
$$
Q(u \vert  A, B, g, k)=A+B\left[1+0.8\times \frac{1-\exp \{-g z(u)\}}{1+\exp \{-g z(u)\}}\right]\left(1+z^{2}(u)\right)^{k} z(u),
$$
where $z(u)=\Phi^{-1}(u)$ is the $u$-th quantile of the standard normal, and $A$, $B>0$, $g$ and $k>-1/2$ are the parameters describing the location, scale, skewness and kurtosis, respectively.
If $g=k=0$, it reduces to a normal distribution. 
ABC is a good candidate for inferring the g-and-k model, as suggested in \citep{AKM:2009}. Here, we apply the adaptive MPMC-IL to approximate the ABC posterior of this model.
Gaussian kernel for ABC is adopted and a mixture of Gaussian distributions is taken as the proposal throughout the example. 

We use the unconstrained parameters ${\theta=(A,\log(B),g,\log(k+1/2))}$ for MPMC-IL and report results for these parameters. The observations in $y_{\text{obs}}$ are independently generated from the g-and-k distribution with fixed parameters $(A,B,g,k)=(3,1,2,0.5)$, from which we get the true unconstrained parameters $\theta_{0}=(3,0,2,0)$.  

In order to obtain complex ABC posterior with multiple modes, the prior density we choose is a Gaussian mixture, i.e., $$p(\theta)=\frac 1 4\sum_{d=1}^{4} N(\theta\vert\mu_d,I_4),$$
$\mu_d=\theta_{0}~+~R_d,$
where $R_d$ is the $d$-th row of a randomly sampled matrix $R$, as shown below:
\begin{equation*}
R=
\left(
\begin{array}{rrrr}
-0.2302 & 0.9273  & 1.3218  & 0.3780   \\
0.0885  & 0.8739  & -0.2305 & -1.0796   \\
-0.8671 & 0.2077  & -0.0338 & 0.4578    \\
0.3725  & -1.0748 & 0.2789  & 0.5326   
\end{array}
\right).\label{eq:summary:octile}
\end{equation*}

\begin{figure}
	\centering    
	\subfigure[Update 1]{\label{ex8:fig:b:1}\includegraphics[width=1.5in]{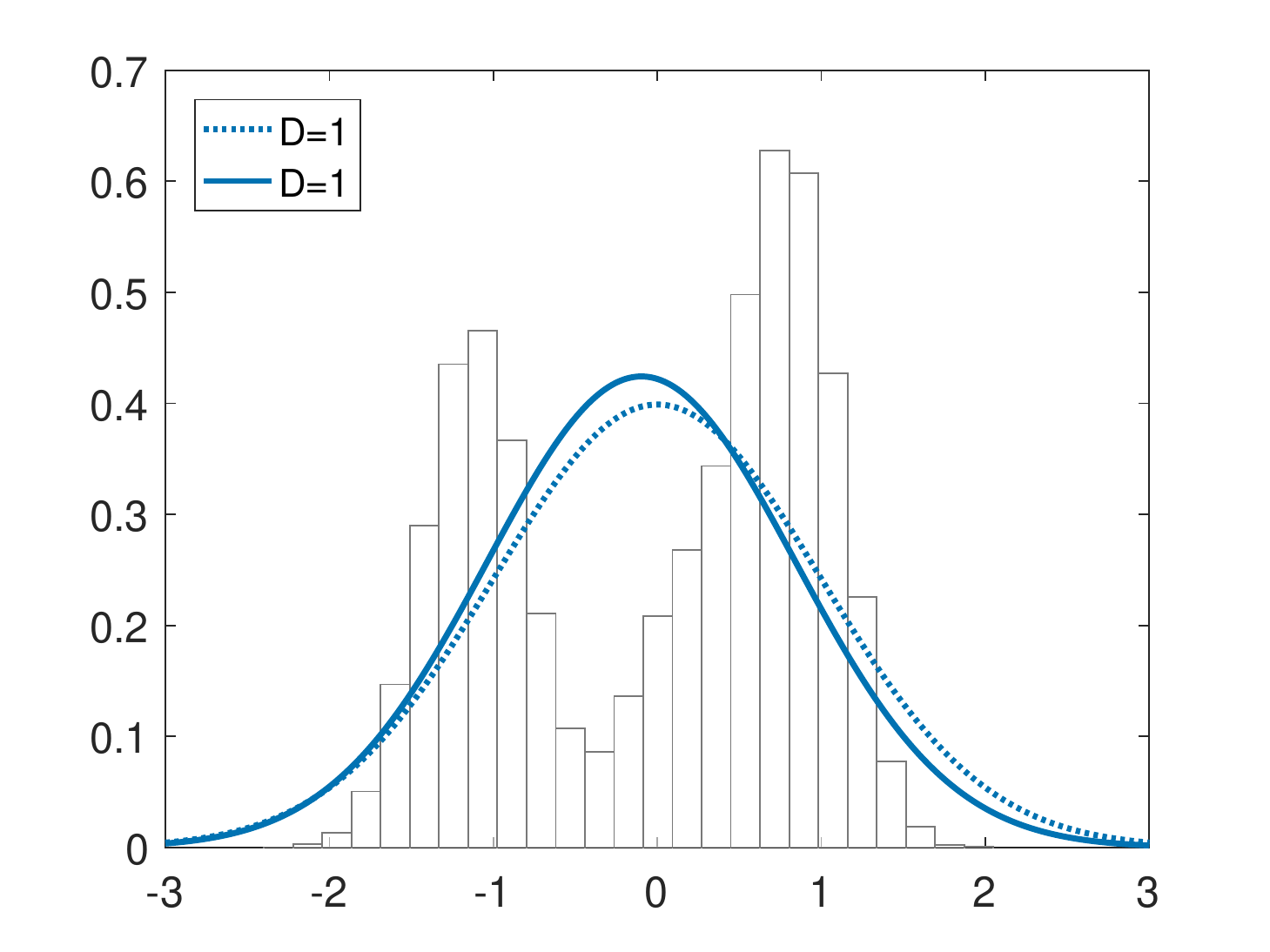}}
	\subfigure[Update 2]{\label{ex8:fig:b:2}\includegraphics[width=1.5in]{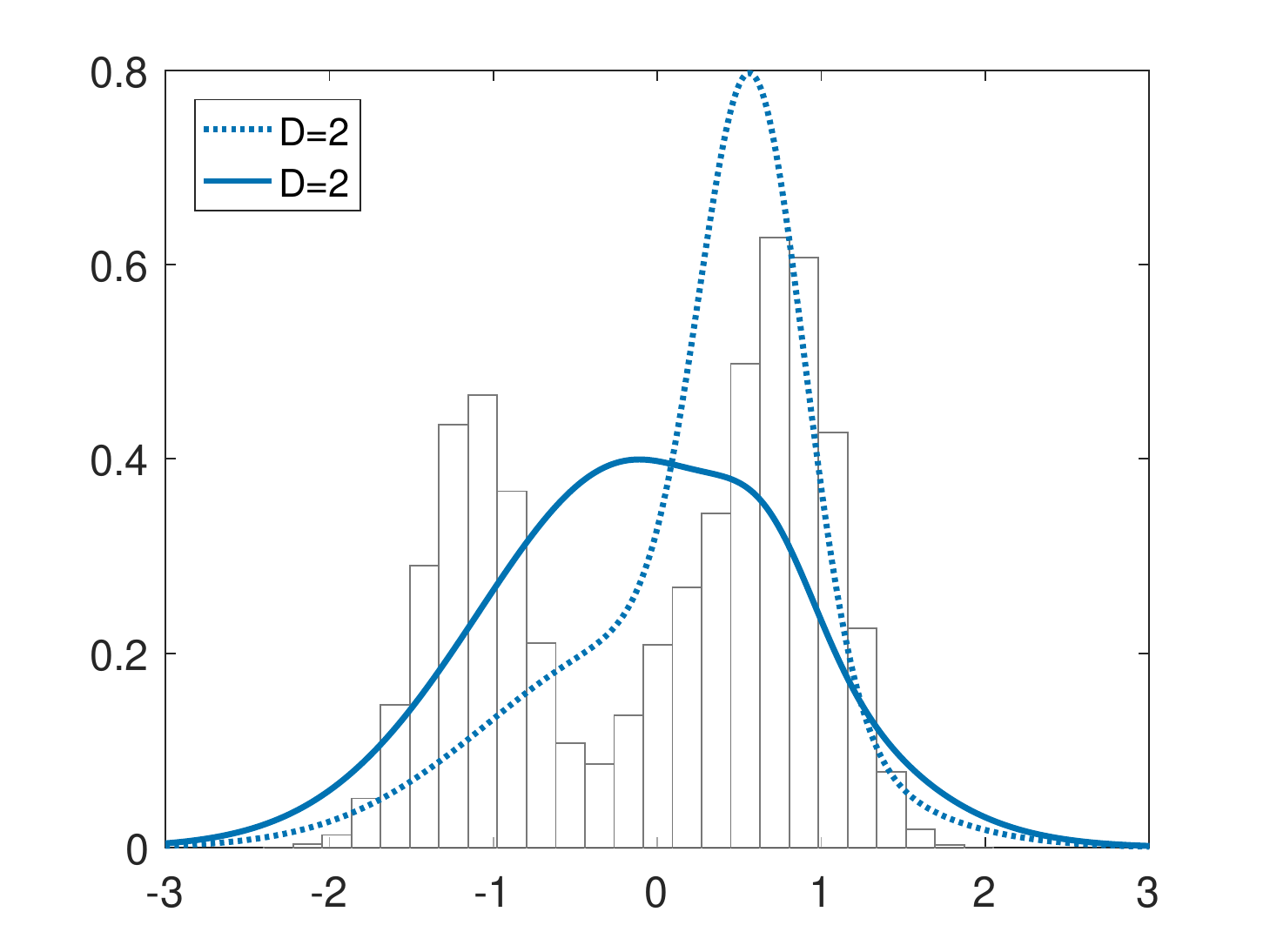}}
	\subfigure[Update 3]{\label{ex8:fig:b:3}\includegraphics[width=1.5in]{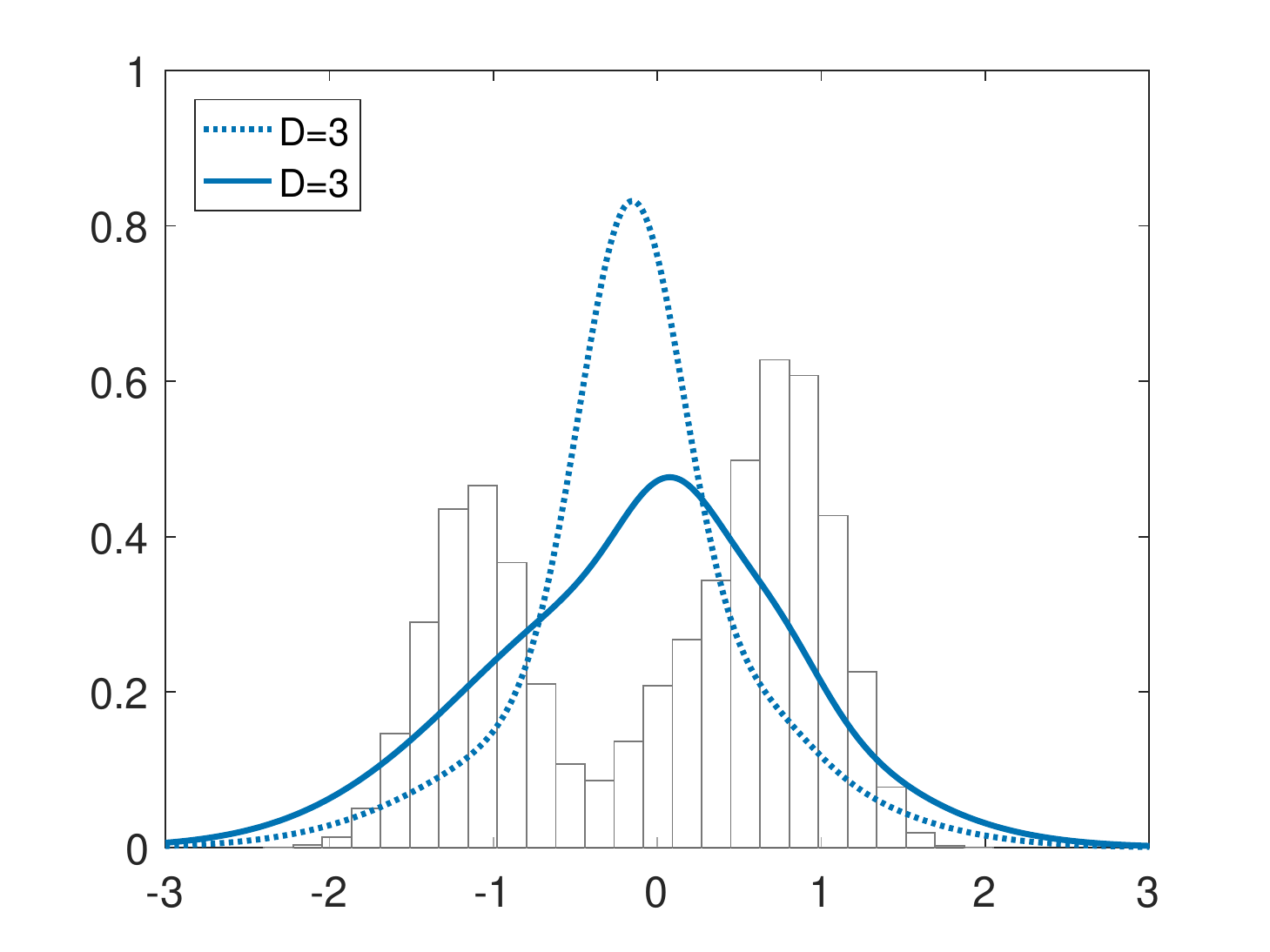}}
	\subfigure[Update 4]{\label{ex8:fig:b:4}\includegraphics[width=1.5in]{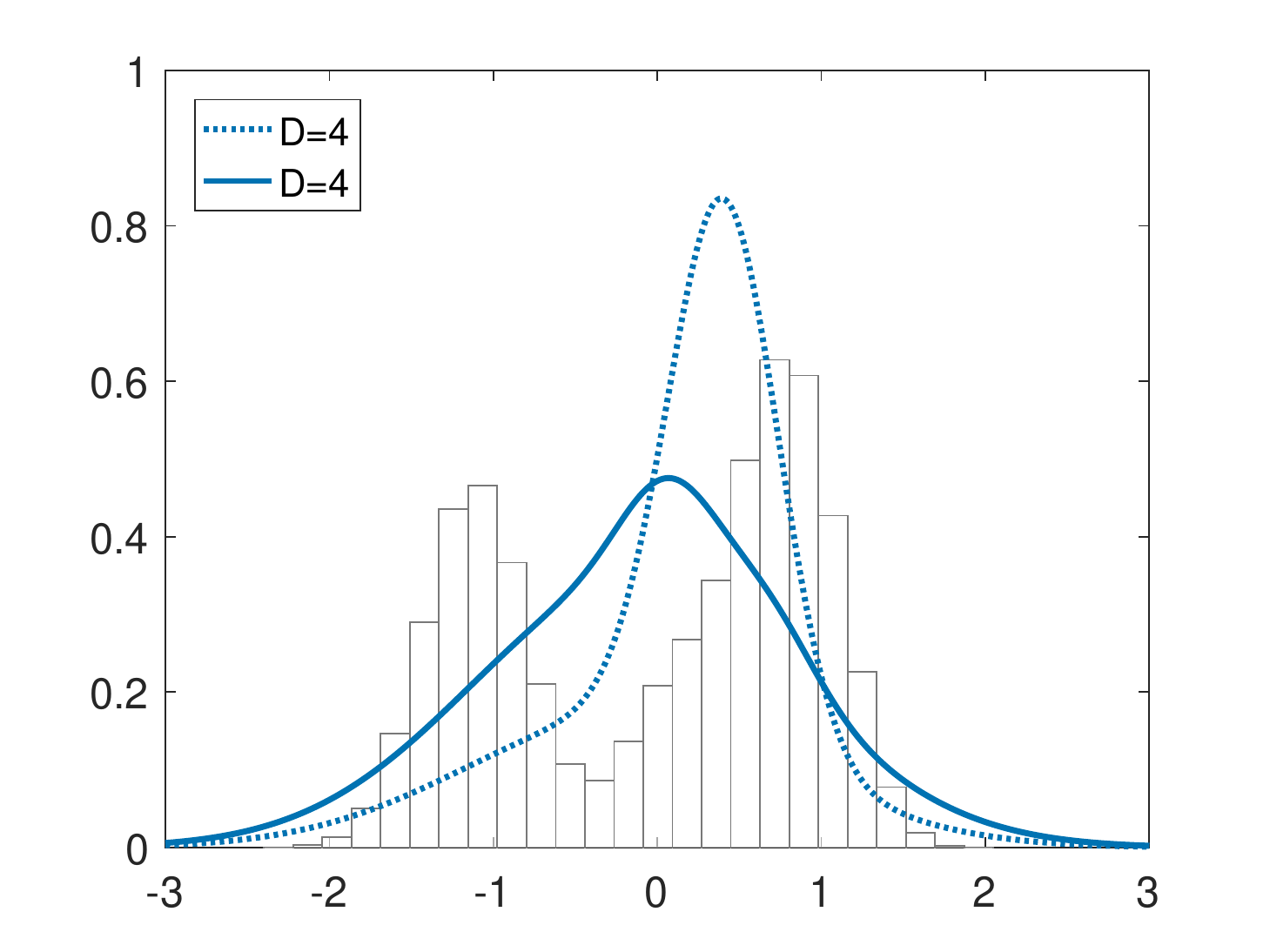}}
	\subfigure[Update 5]{\label{ex8:fig:b:4}\includegraphics[width=1.5in]{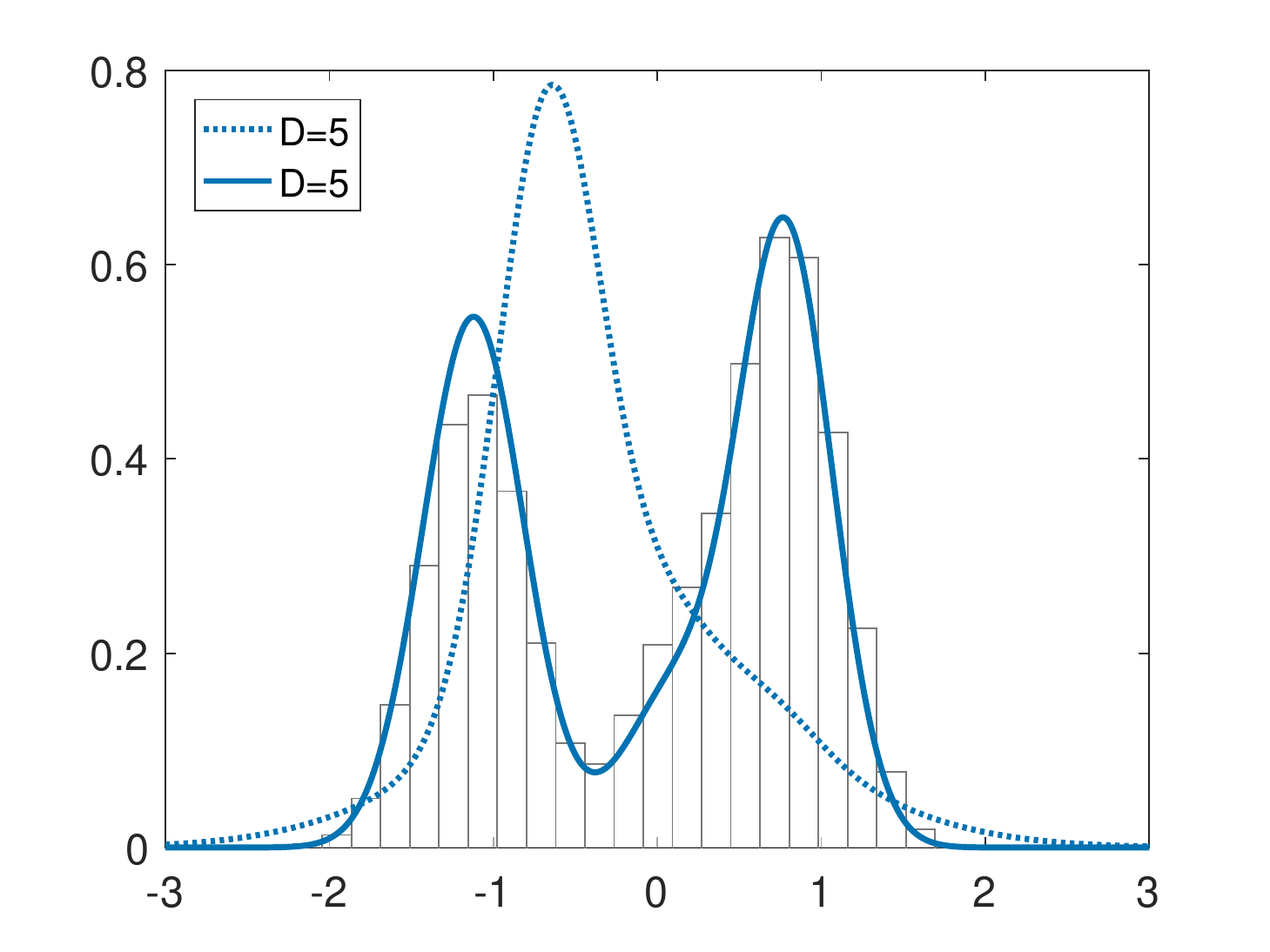}}
	\subfigure[Update 6]{\label{ex8:fig:b:4}\includegraphics[width=1.5in]{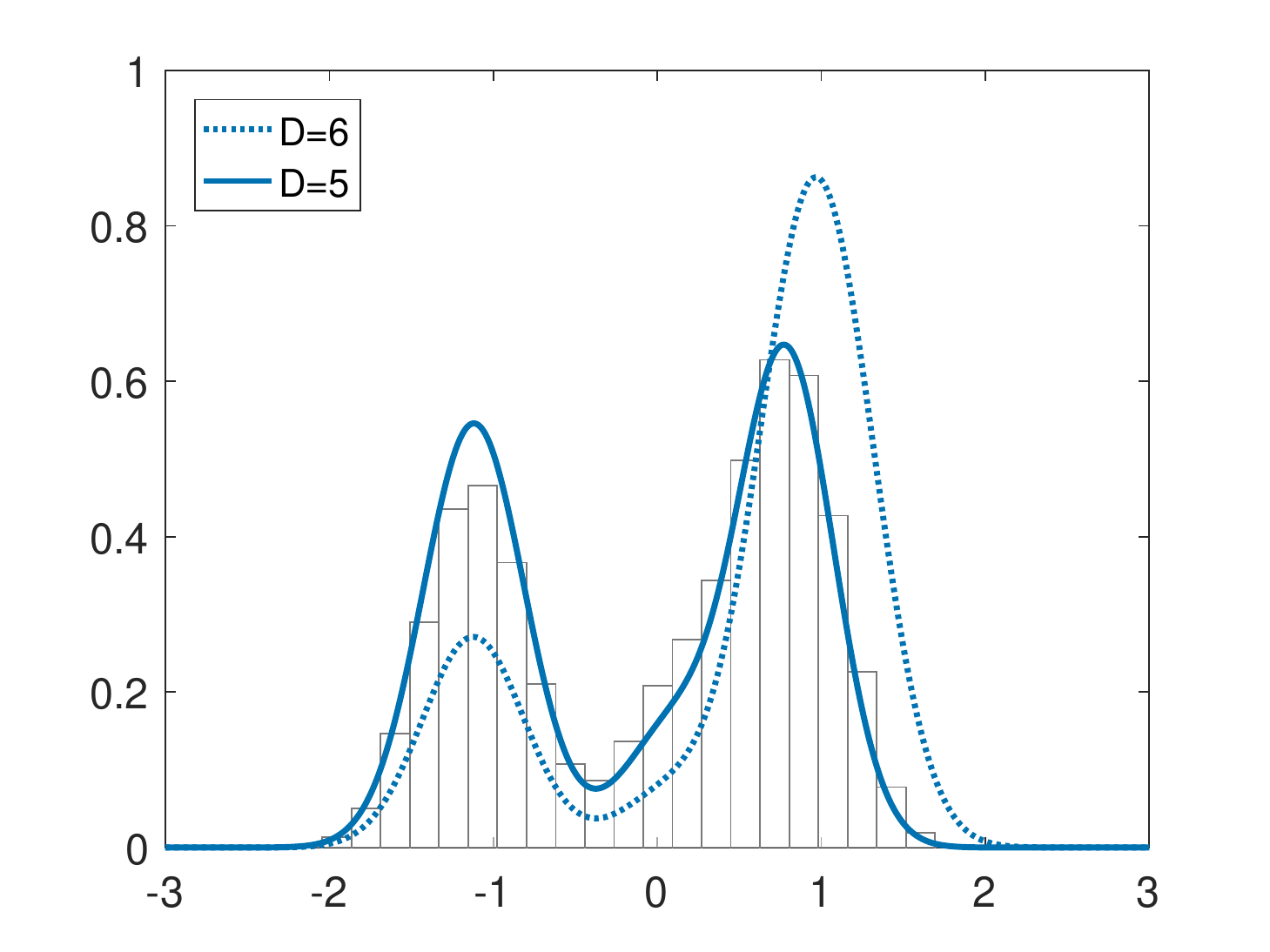}}
	\caption{Approximate marginal posterior distributions of the parameter $\theta_2 = \log(B)$ for six successive updates of MPMC-IL with  a fixed window rule $t_w=20$. The histograms are obtained from ABC acceptance-rejection samples of length $10^5$. The dashed lines are the initial mixtures, and the solid lines are the outputs of each update. The settings of ABC in this test: $\mathrm{dim}(\yobs)=20$, $h=12.34$ and $S(y)=y$.		
	}\label{fig:ex8:abgk:b:20}
\end{figure}

First, we look at the effect of a small data set on our proposed algorithm in which the number of observations is $20$. For this case, we take the entire data set as the summary statistics (i.e., $S(y)=y$) so that the ABC posterior converges to the exact posterior as the bandwidth $h\to 0$. We set $h=12.34$ in ABC.
The initial distribution for the adaptive MPMC-IL algorithm is a standard normal with $D=1$. We set $N=N_{\text{add}}=10^5$, and a fixed window size $t_w=20$ in Algorithm~\ref{algo:component--updating}. Figure \ref{fig:ex8:abgk:b:20} illustrates the marginal distribution of the parameter $\theta_2 = \log(B)$ derived from the adaptive MPMC-IL algorithm with six successive updates in the component--updating procedure. The initial and the resulting mixture distributions  for each update are depicted, along with the ABC acceptance-rejection (the benchmark) with $10^5$ samples. 

From the comparison of the approximate posteriors in different updates, we see an evident improvement of the posterior by the MPMC-IL algorithm, where a distribution closer to the benchmark is observed. Moreover, with the increasing number of components in Gaussian mixture, the derived mixture distribution is progressively improved with finer details. It is worth noting that at the last update in Figure \ref{fig:ex8:abgk:b:20}-(f), the derived mixture has $5$ components other than $6$, because a component with small weight was removed after performing the main MPMC-IL algorithm. Such adaptive handling of components in the proposal helps to avoid over-fitting and unnecessary calculations in Gaussian mixture, thus improving the stableness and efficiency of the proposed algorithm. 

\begin{figure}[htb]
	\centering	\includegraphics[width=0.5\linewidth]{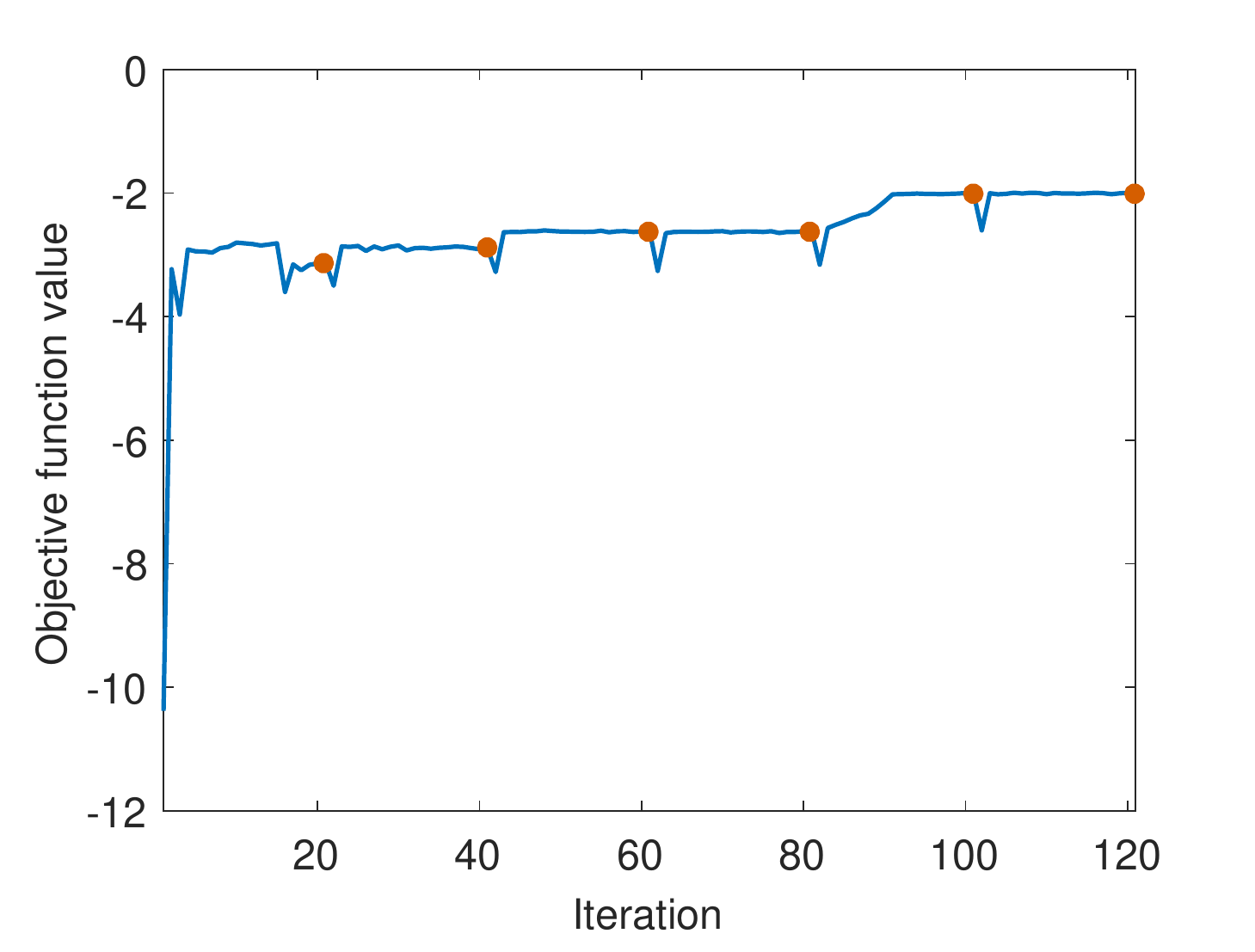}
	\caption{The objective function values over iterations for the adaptive MPMC-IL algorithm. The component number is successively updated after every fixed window size ($t_w = 20$) of iterations. The settings are the same as in Figure~\ref{fig:ex8:abgk:b:20}.  The highlighted points (red) indicate the iterations at which MPMC-IL stops for each update. 
	}	\label{fig:ex8:lb:obs20:summary1}
\end{figure}

Figure \ref{fig:ex8:lb:obs20:summary1} monitors the convergence of the adaptive MPMC-IL algorithm with the fixed window ($t_w=20$) stopping rule \eqref{eq:fixedwindow}. 
More stable convergences of the objective function values are observed with larger $D$. Moreover, from the overall look of the objective function values in Figure \ref{fig:ex8:lb:obs20:summary1}, we can see progressive improvements of the objective function values after successive updates of the number of components $D$, which indicates the effectiveness of the proposed component--updating procedure. It should be noted that there is an evident sharp decrease of the objective function after each updating of the number $D$. The reason is that a newly added component would usually destroy the current good proposal derived from the MPMC-IL algorithm and might leads to a lower value. Since the new mixture component is initialized adaptively, the objective function value does not decline much. This highlights the importance of initialing a  new mixture component. A quick increase and then a nearly stable trend of the objective function values within each update after the first iteration are observed. Such fast approaching for the state of convergence shows the effectiveness of the MPMC-IL algorithm.

Next, we look at the effect of a large data set on our proposed algorithm in which the number of observations is increased from $20$ to $1000$. For this case, following \citep{drovandi:2011}, we take a low-dimensional summary statistics given by
${S(y)=\left(S_{A}, S_{B}, S_{g}, S_{k}\right)}$ with
\begin{equation}\label{quantile:ss}
\begin{aligned}
S_{A}&=E_{4},\\
S_{B}&=E_{6}-E_{2},\\
S_{g}&=\left(E_{6}+E_{2}-2 E_{4}\right)/S_{B}, \\ 
S_{k}&=\left(E_{7}-E_{5}+E_{3}-E_{1}\right) / S_{B},
\end{aligned}
\end{equation}
where
${E_{1} \leq E_{2} \leq \cdots \leq E_{7}}$ are the octiles of ${y}$.
We set the bandwidth $h = 0.5971$.

\begin{figure}
	\centering    
	\subfigure[$\theta_1=A$]{\label{ex7:fig:a}\includegraphics[width=2in]{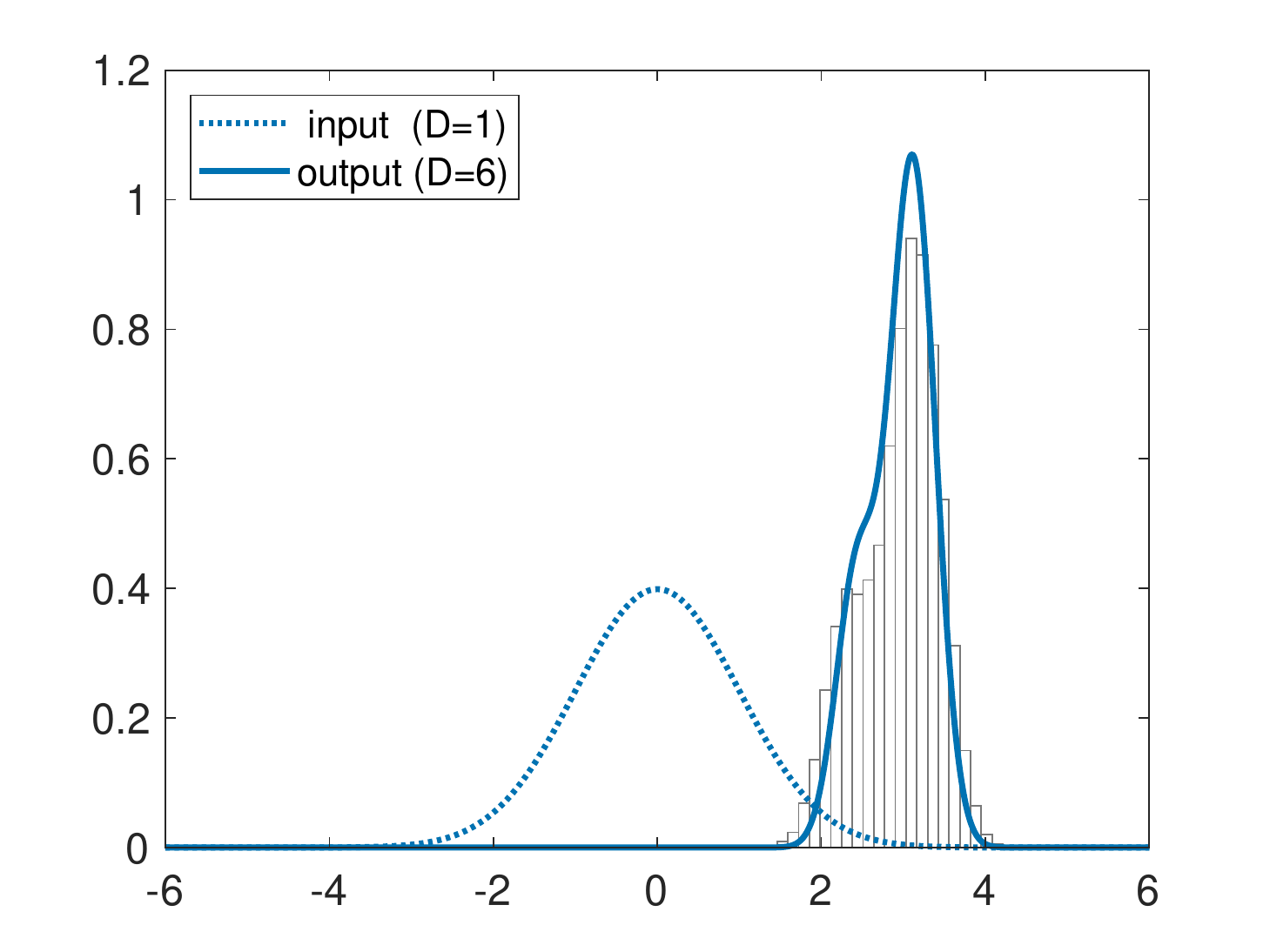}}
	\subfigure[$\theta_2=\log(B)$]{\label{ex7:fig:b}\includegraphics[width=2in]{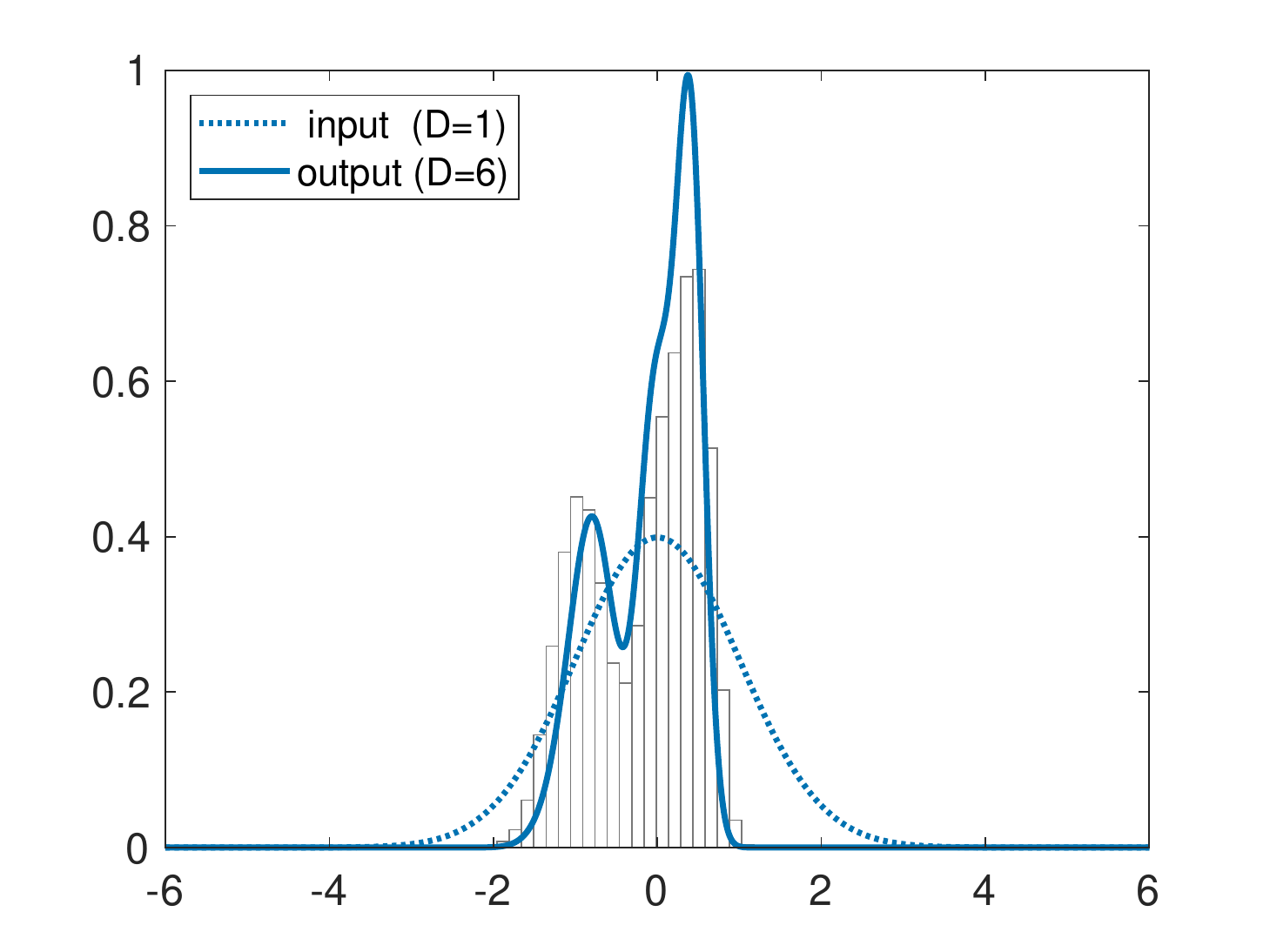}}
		\subfigure[$\theta_3=g$]{\label{ex7:fig:g}\includegraphics[width=2in]{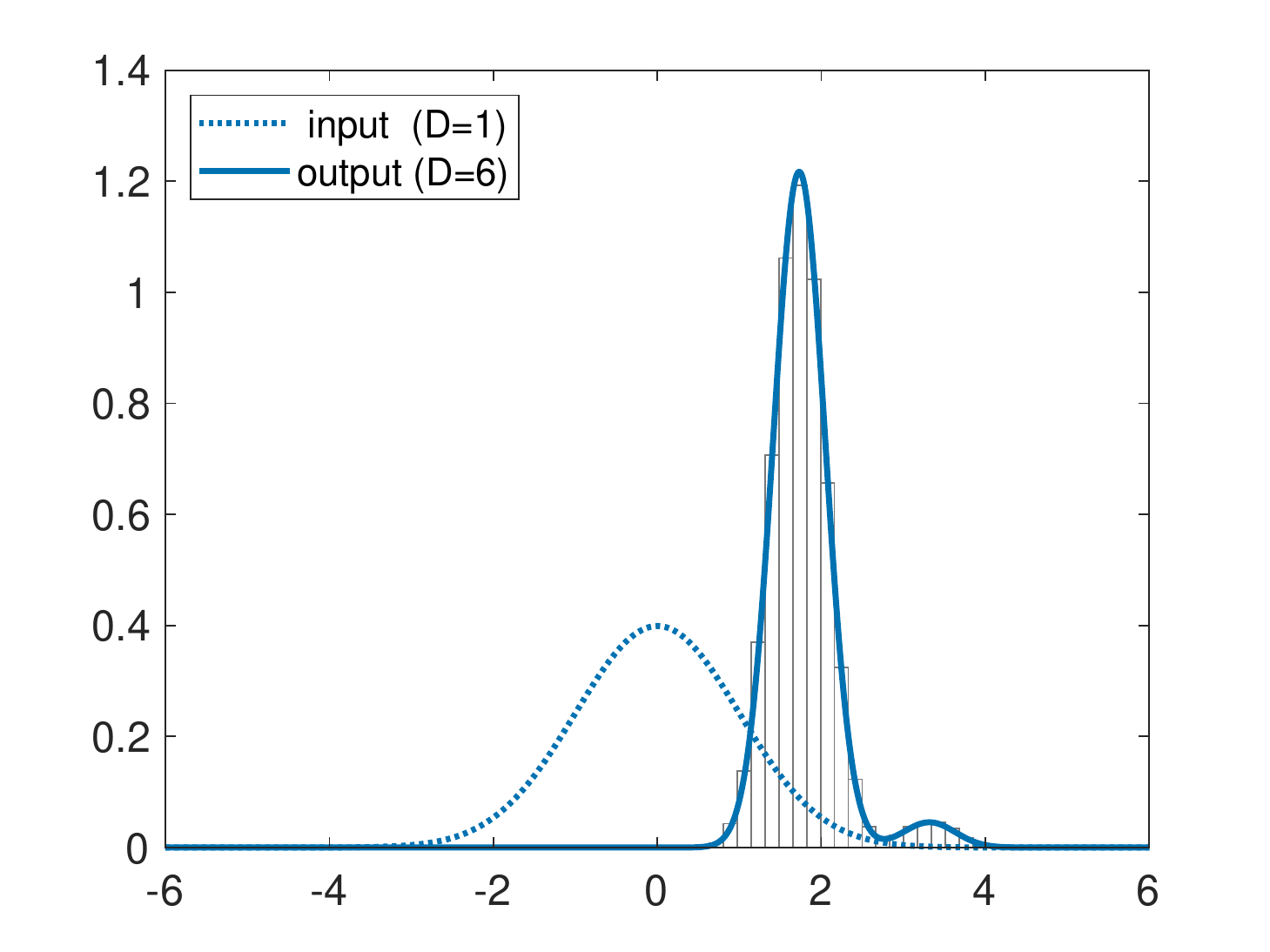}}
	\subfigure[$\theta_4=\log(k+1/2)$ ]{\label{ex7:fig:k}\includegraphics[width=2in]{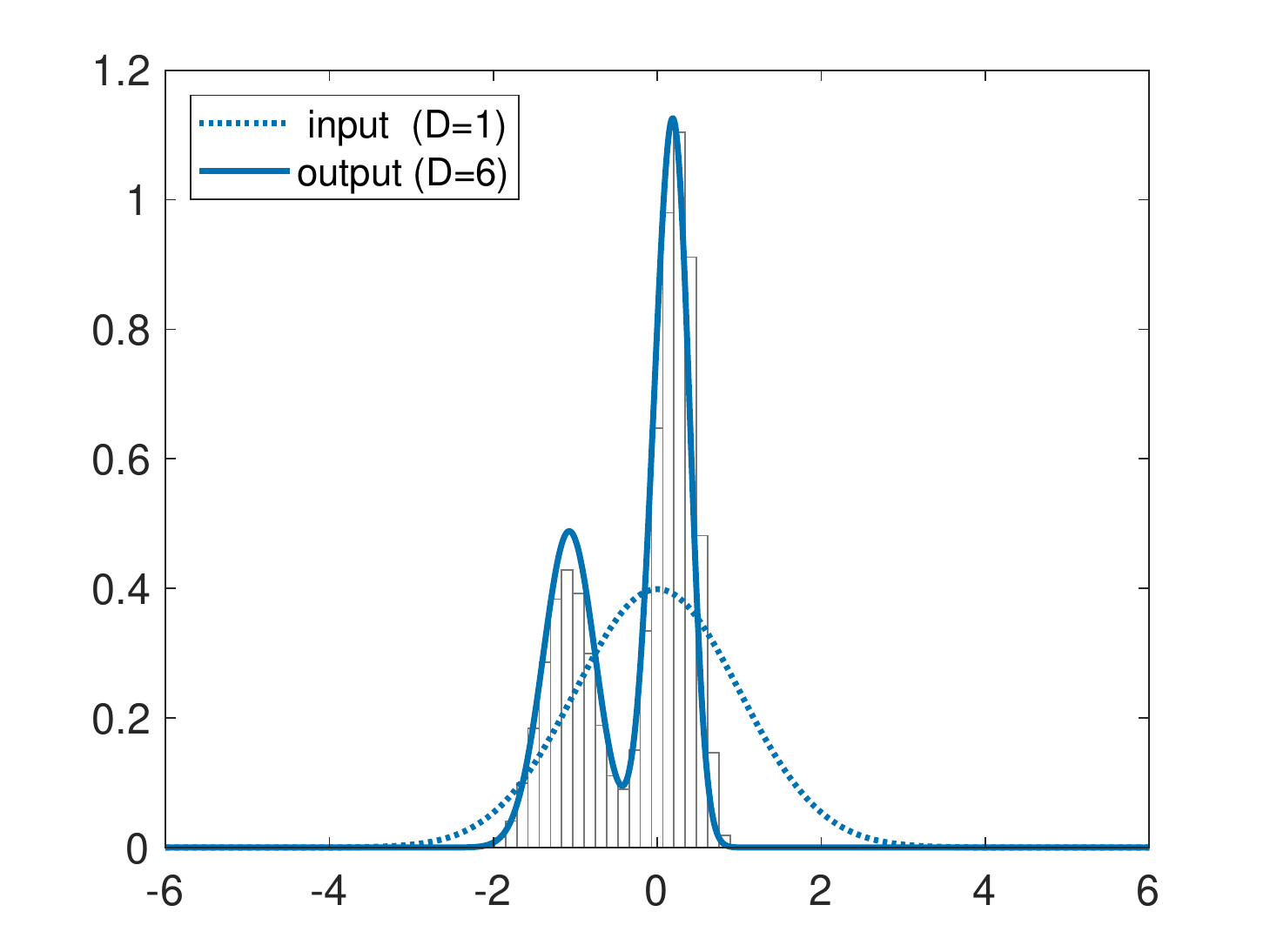}}
	\caption{Approximate marginal posterior distributions of $A$, $\log(B)$, $g$ and $\log(k+1/2)$ derived from the adaptive MPMC-IL algorithm after six updates with a fixed window rule $t_w=20$. The histograms are obtained from ABC acceptance-rejection samples of length $10^5$. The dashed lines are the initial distributions (standard Gaussian), and the solid lines are the outputs of MPMC-IL after six updates ($D=6$). The settings of ABC in this test: $\mathrm{dim}(\yobs)=1000$, $h=0.5971$ and $S(y)$ is given by \eqref{quantile:ss}.		
		}\label{fig:ex7:abgk:1e3}
\end{figure}

\begin{figure}[htb]
	\centering
	\includegraphics[width=0.6\linewidth]{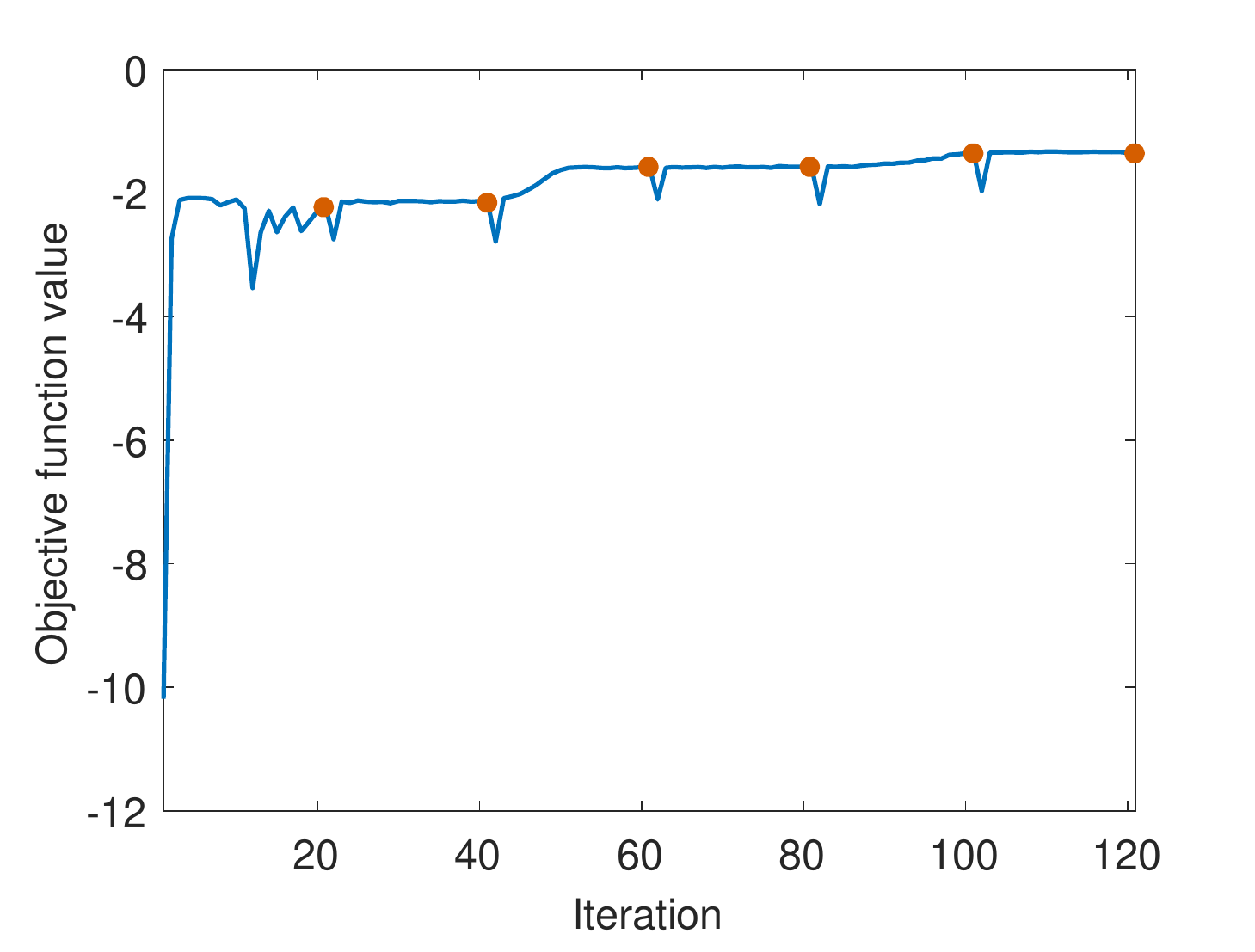}
	\caption{The objective function values over iterations for the adaptive MPMC-IL algorithm. The component number is successively updated after every fixed window size ($t_w = 20$) of iterations. The settings are the same as in Figure~\ref{fig:ex7:abgk:1e3}.  The highlighted points (red) indicate the iterations at which MPMC-IL stops for each update.
	}\label{fig:ex7:lb:obs1000:summary2}
\end{figure}

Figure \ref{fig:ex7:abgk:1e3} illustrates the approximate marginal posterior distributions for the four unconstrained parameters ${\theta=(A,\log(B),g,\log(k+1/2))}$ resulting from the adaptive MPMC-IL algorithm after six updates of $D$ with a fixed window rule $t_w=20$, respectively. The objective function values are provided in Figure~\ref{fig:ex7:lb:obs1000:summary2}. One sees better approximation for each parameter compared with the  initial distributions (i.e., standard Gaussian) in Figure \ref{fig:ex7:abgk:1e3}, and the resulting mixture distributions fit the benchmarks very well, exploiting the multimodality of the ABC posterior. In Figure \ref{fig:ex7:lb:obs1000:summary2}, we observe a similar pattern as in Figure~\ref{fig:ex8:lb:obs20:summary1}. Choosing a proper $D$ is critical for MPMC-IL, and our adaptive procedure in Algorithm~\ref{algo:component--updating} accomplished this goal with a large value of the objective function. We should note that the adaptive MPMC-IL algorithm starts from a ``bad" initial distribution for $D=1$, which is far way from the benchmark distribution. The MPMC-IL algorithm then converges quickly due to the goodness of the initialized mixture distribution for each round of update.

%

\begin{figure}[htb]
	\centering	\includegraphics[width=.6\linewidth]{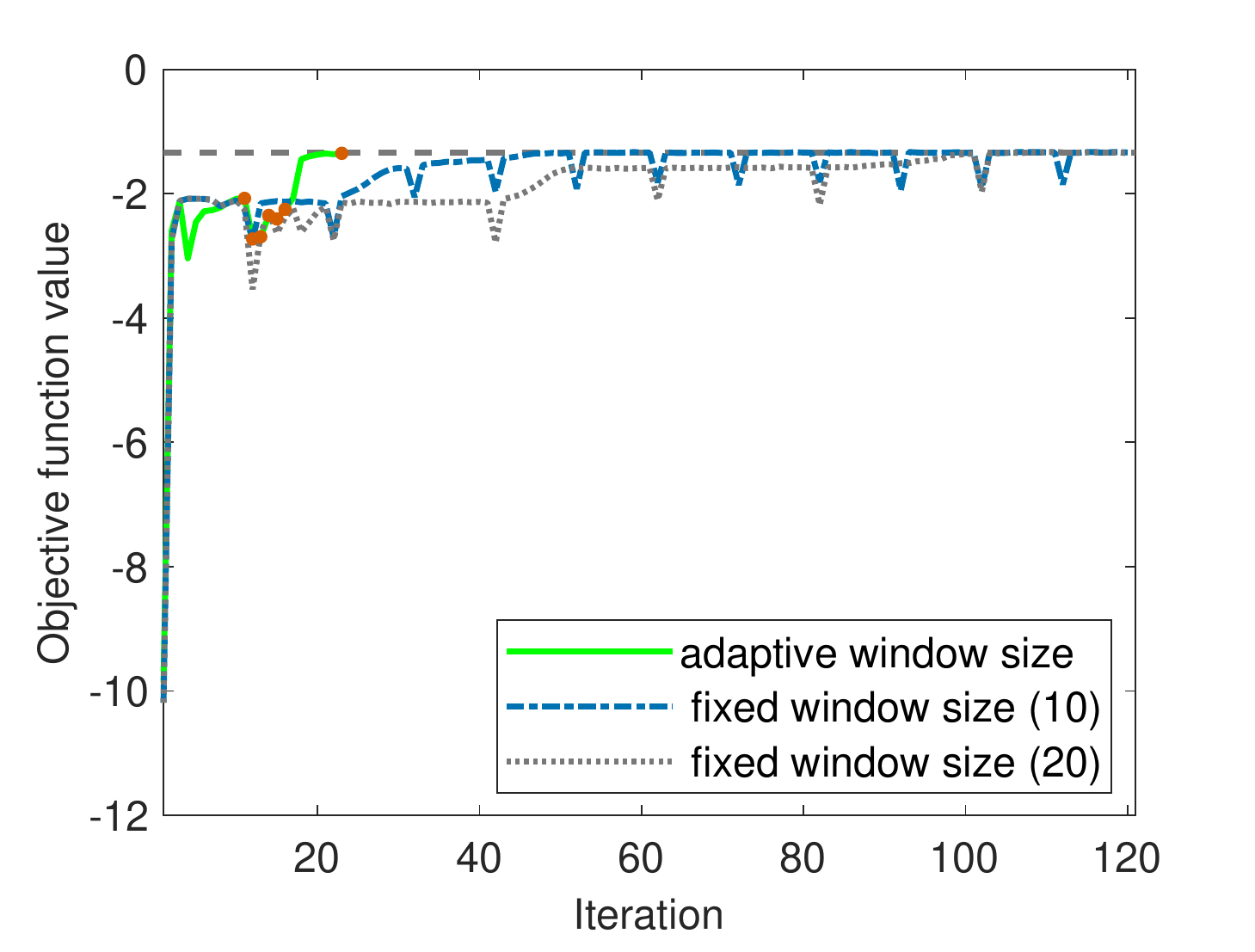}
	\caption {Comparison of the objective function values versus the iterations with the fixed window size rule \eqref{eq:fixedwindow} and the adaptive window size rule \eqref{eq:adaptwindow}. The highlighted points (red) indicate the iterations at which MPMC-IL stops under the adaptive paradigm. The settings are the same as in Figure~\ref{fig:ex7:abgk:1e3}. The fixed window sizes are $t_w=10$ and $t_w=20$, respectively. The parameters used for the adaptive window size rule \eqref{eq:adaptwindow} are $s=5$, $\epsilon_0=0.1$.}\label{fig:1000pts:windowsize:overall}
\end{figure}

In the previous experiments, we used the fixed window size rule \eqref{eq:fixedwindow} to update the number of component $D$. Now we examine the effect of the window size $t_w$ and compare the fixed window size rule \eqref{eq:fixedwindow} with the adaptive window size rule \eqref{eq:adaptwindow} for the g-and-k model with $1000$ observations. 
Figure \ref{fig:1000pts:windowsize:overall} depicts the comparison of the objective function values over the iterations with the two updating rules, where $t_w=10$ and $20$, respectively.
We can see from the figure that the adaptive paradigm drives more frequent updates than the fixed $t_w$ paradigm, hence achieving the largest objective function values with the fewest iterations.
For this example, MPMC-IL with the smaller fixed window size $t_w=10$ converges faster than the other case of fixed window size.
This indicates that $t_w$ has an impact on the convergence of the proposed MPMC-IL algorithm, and the adaptive window size rule \eqref{eq:adaptwindow} indeed helps to speed up the algorithm.

%
%

\subsection{Generalized linear mixed model}

Generalized linear mixed models (GLMM) use random effects $\alpha_i$ to account for the dependence between the observations $y_i=\{y_{ij},~j=1,\dots,n_i\}$, which are measured on the same individual $i=1,\dots,n$. Now $\yobs=\{y_1,\dots,y_n\}$. The joint distribution of the observable $y_i$ and the (unobservable) random effect $\alpha_i$ is given by $$p(y_i,\alpha_i\vert\theta)= p(\alpha_i\vert\theta)p(y_i\vert \alpha_i,\theta),$$ 
where $\theta$ is the model parameter of interest. In practice, the conditional distributions $p(\alpha_i\vert\theta)$ and $p(y_i\vert\theta,\alpha_i)$ are tractable, but  
$$p(y_i\vert\theta)=\int p(y_i\vert\alpha_i,\theta)p(\alpha_i\vert\theta)\mrd \alpha_i$$
is analytically intractable in most cases. This makes the likelihood of the $n$ independent observations $p(\yobs\vert\theta)=\prod_{i=1}^n p(y_i\vert\theta)$ intractable. Suppose that we are able to generate samples from the distributions $p(\alpha_i\vert\theta)$. The likelihood $p(\yobs\vert\theta)$ is estimated unbiasedly by
\begin{equation}
\hat p(\yobs\vert\theta)=\prod_{i=1}^n \hat p_{N_i}(y_i\vert\theta)
\end{equation}
with
\begin{equation}\label{eq:liki}
\hat p_{N_i}(y_i\vert\theta)=\frac{1}{N_i}
\sum_{j=1}^{N_i}p(y_i\vert\alpha_i^{(j)},\theta),
\end{equation}
where $\alpha_i^{(j)}\simiid p(\alpha_i\vert\theta)$, $j=1,\dots,N_i$. In our setting, $\Psi(X;\yobs,\theta)=\hat p (\yobs\vert\theta)$ so that $\mbe[\Psi(X;\yobs,\theta)\vert\theta]=p(\yobs\vert\theta)$.

We consider the logistic regression with Six City data set \cite{Fitz:1993}, which was studied in \cite{tran:2017} and \cite{he:2021} in the context of variational Bayes. The data consist of binary responses $y_{ij}$ which is the wheezing status (1 if wheezing, 0 if not wheezing) of the $i$th child at time-point $j$, where $i=1,\dots,537$ (representing 537 children) and $j=1,2,3,4$ (corresponding to $7,8,9,10$ year-old). Covariates are $A_{ij},$ the age of the $i$th child at time-point $j$, and $S_i$ the $i$th maternal smoking status (0 or 1). For the numerical stability, the covariates $A_{ij}$ are centered at 9 years, i.e., $A_{ij}=j-3$ for $j=1,\dots,4$ and all $i$. We consider the logistic regression model with a random intercept: $y_{ij}\vert\beta,\alpha\sim \mathrm{Bernoulli}(p_{ij})$, where $$p_{ij}=\frac{\exp(\beta_1+\beta_2A_{ij}+\beta_3S_i+\alpha_i)}{1+\exp(\beta_1+\beta_2A_{ij}+\beta_3S_i+\alpha_i)}$$ and $\alpha_i\sim N(0,\tau^2)$. The  model parameters are $\beta_1,\beta_2,\beta_3$, and $\tau^2$. Note that
$$p(y_i\vert\alpha_i,\theta) =\prod_{j=1}^4\frac{\exp\{y_{ij}(\beta_1+\beta_2A_{ij}+\beta_3S_i+\alpha_i)\}}{1+\exp(\beta_1+\beta_2A_{ij}+\beta_3S_i+\alpha_i)}.$$ 
So it suffices to generate $\alpha_i^{(j)}\sim N(0,\tau^2)$ independently in \eqref{eq:liki}.

Again, we use unconstrained parameters $ \theta=(\beta_1,\beta_2,\beta_3,\log\tau^2) $ for our proposed algorithm, and take independent priors: $\beta_i\sim N(0,50)$ and $\tau\sim \mathrm{Gamma}(1,0.1)$. The initial distribution for the adaptive MPMC-IL algorithm is also a standard Gaussian ($D=1$). We set $N=N_{\text{add}}=10^4$, and a fixed window size $t_w=10$ in Algorithm~\ref{algo:component--updating}. To estimate the likelihood, we use $N_i=500$ in the estimator \eqref{eq:liki}  for all $i$. The RStan package `rstanarm' is used to generate MCMC samples from the posterior $p(\theta\vert\yobs)$ as a benchmark, which performs posterior analysis for models with dependent data such as GLMMs.

\begin{figure}
	\centering
	\includegraphics[scale=0.5]{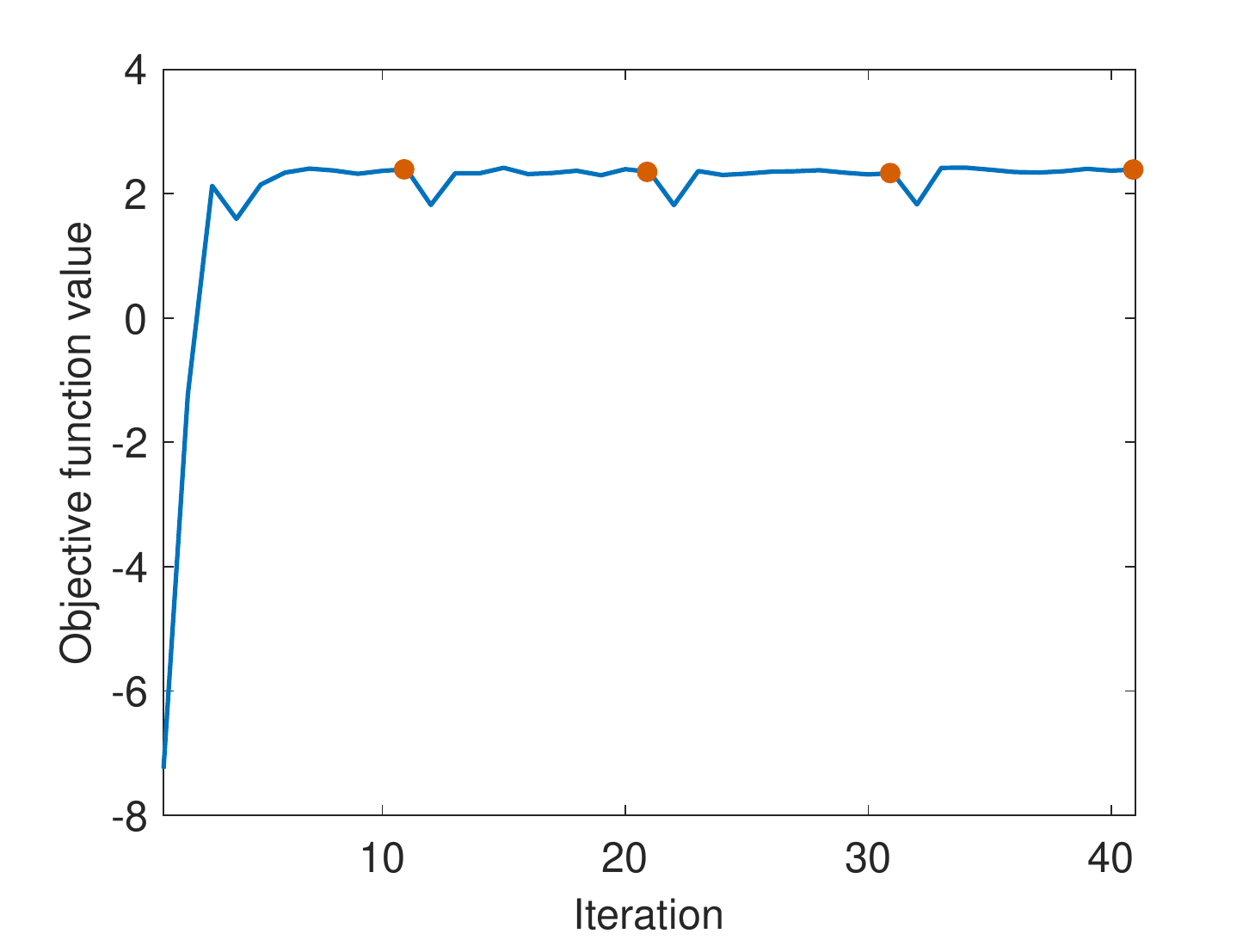}
	\caption{The objective function values over iterations for the GLMM model. The component number is successively updated after every fixed window size ($t_w = 10$) of iterations. The highlighted points (red) indicate the iterations at which MPMC-IL stops for each update.}	
	\label{LB}
\end{figure}

\begin{figure}
	\centering    
	\subfigure[$\theta_1=\beta_1$]{\label{sixcity:beta1}\includegraphics[width=2in]{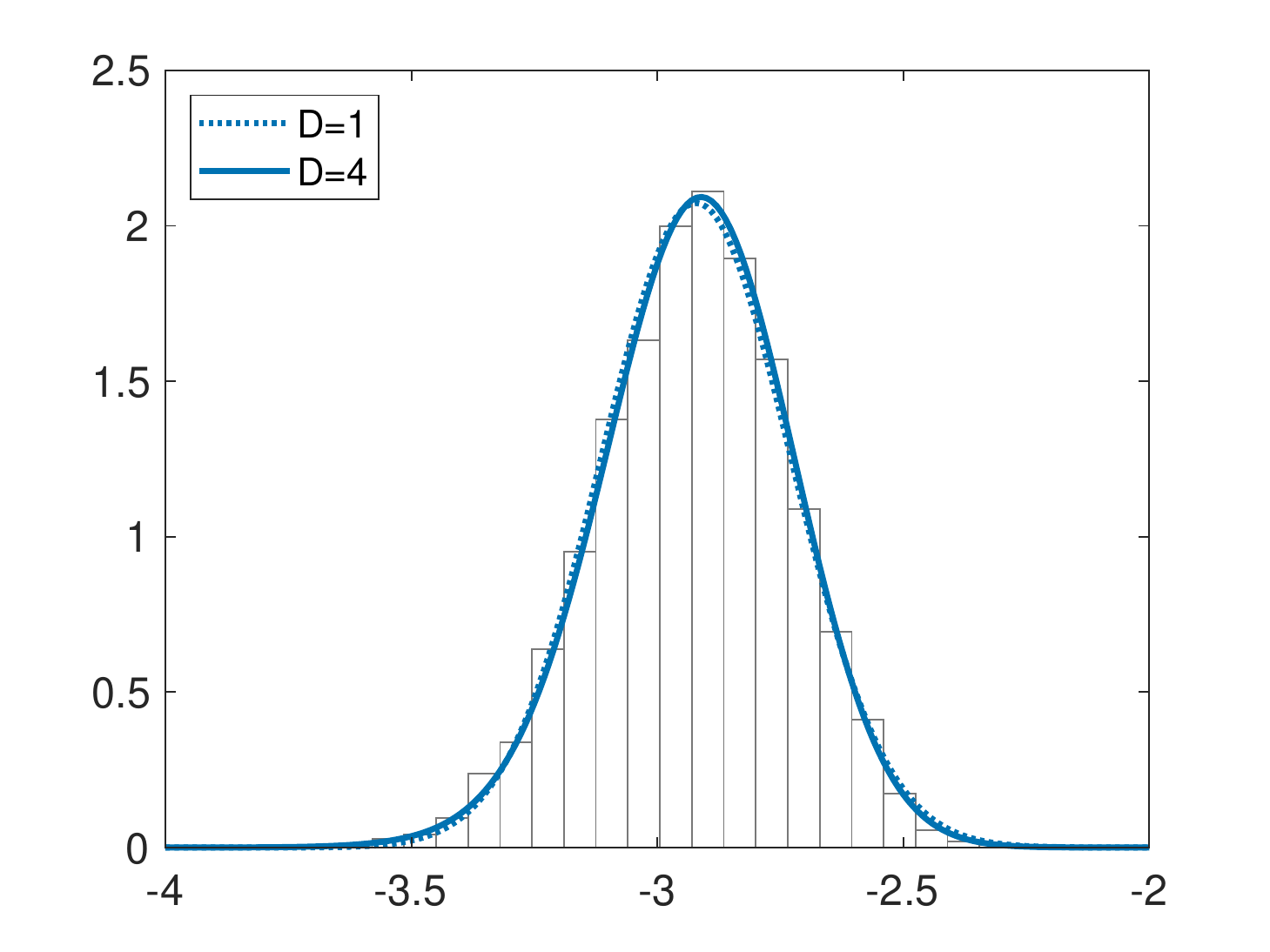}}
	\subfigure[$\theta_2=\beta_2$]{\label{sixcity:beta2}\includegraphics[width=2in]{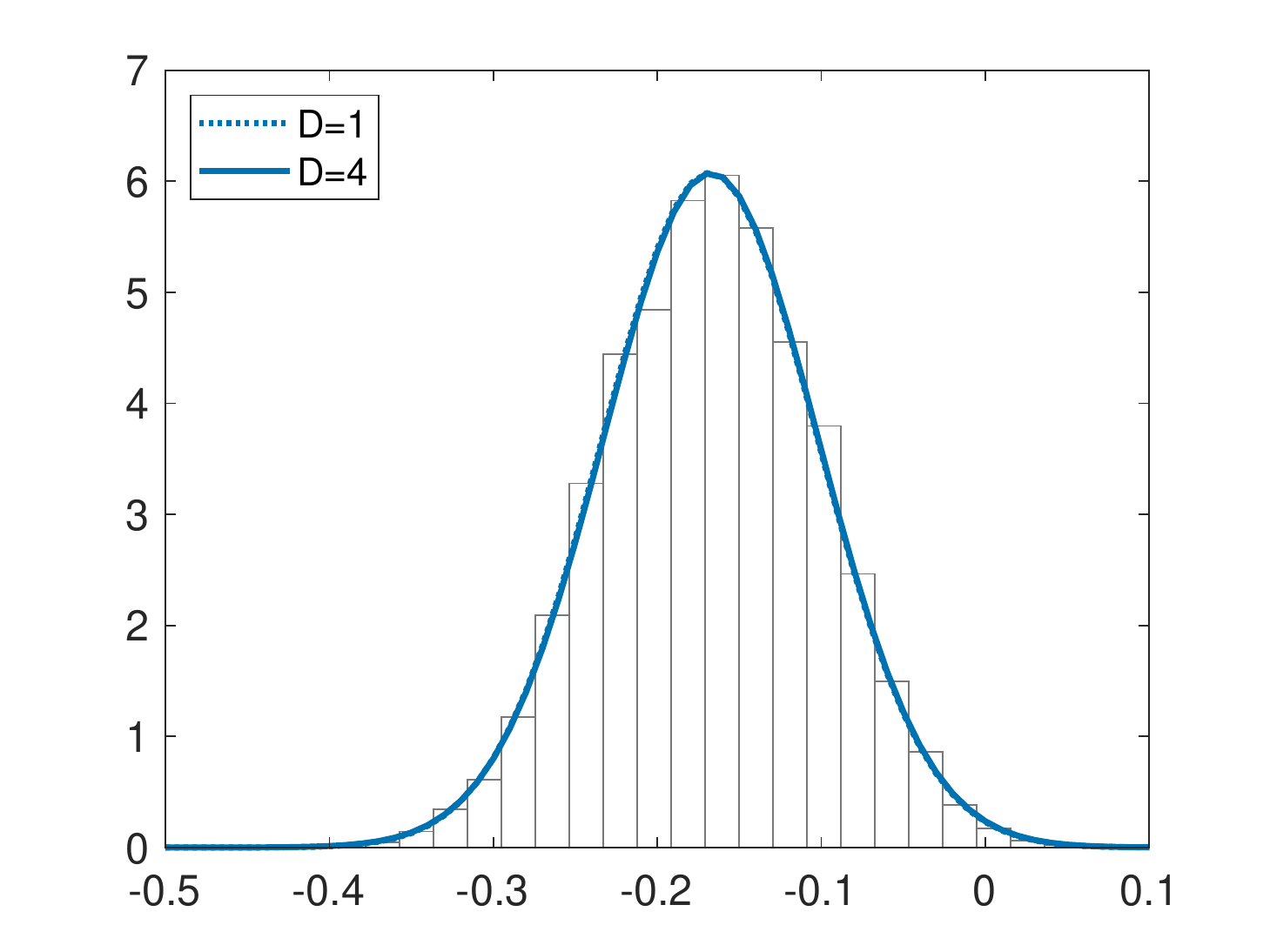}}
	\subfigure[$\theta_3=\beta_3$]{\label{sixcity:beta3}\includegraphics[width=2in]{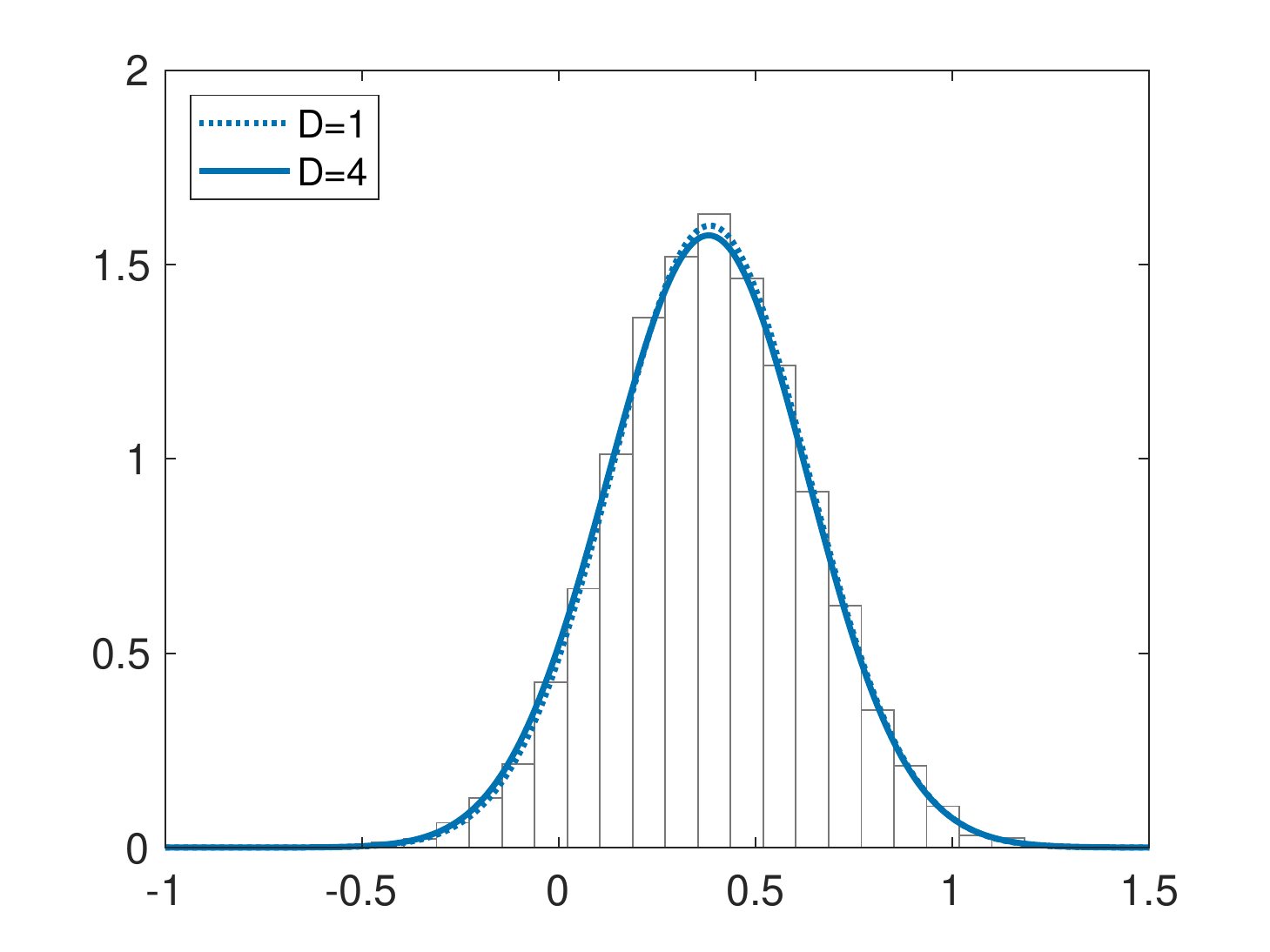}}
	\subfigure[$\theta_4=\log \tau^2$ ]{\label{sixcity:logtau2}\includegraphics[width=2in]{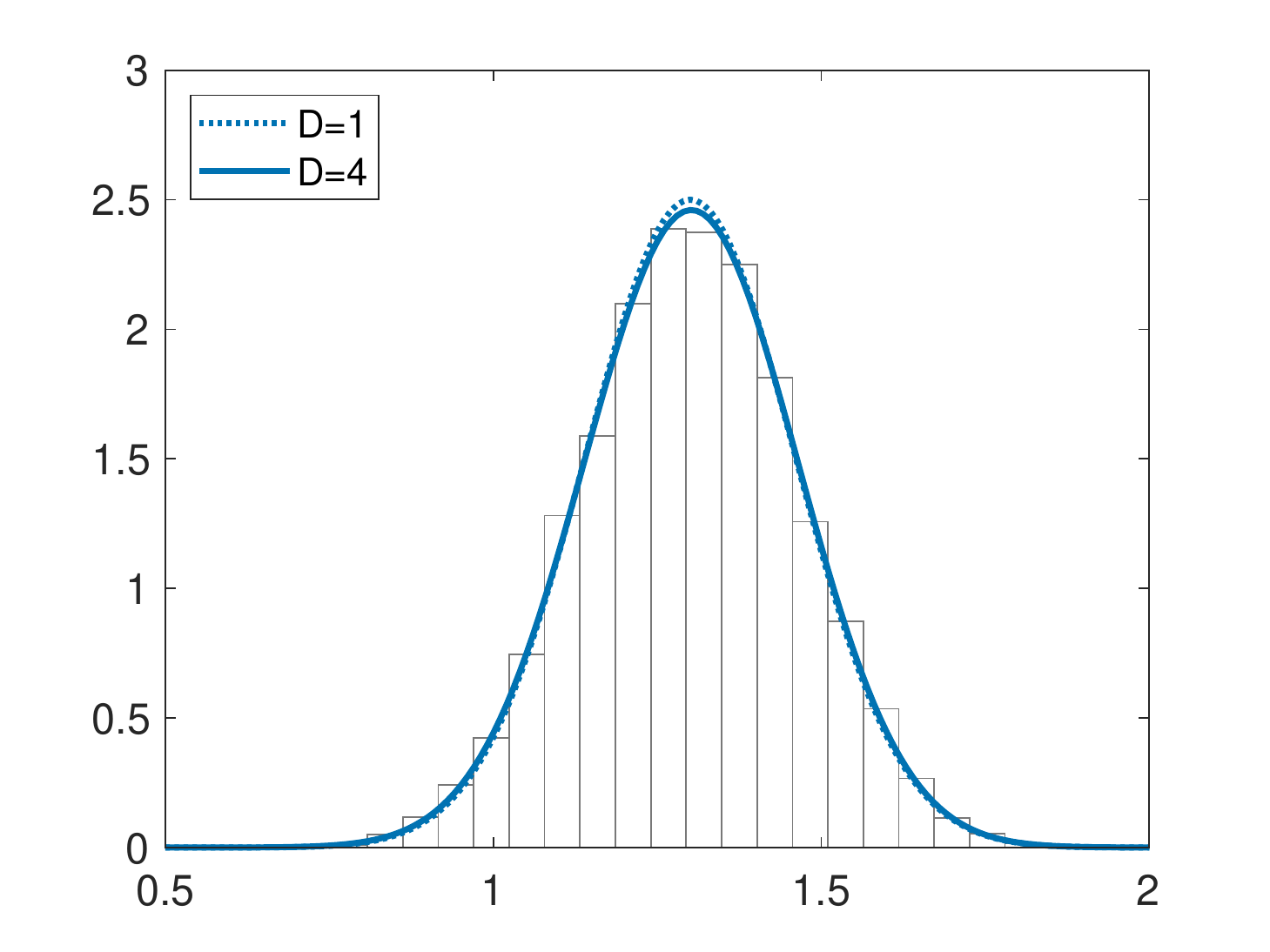}}
	\caption{Comparison of approximate marginal posterior distributions of $ (\beta_1,\beta_2,\beta_3,\log \tau^2) $ derived from the adaptive MPMC-IL algorithm with update $1$ and update $4$ under the fixed window rule $t_w=10$. The histograms are obtained from MCMC samples of length $10^5$ by using RStan package `rstanarm'. The dashed lines and the solid lines are the outputs of $D=1$ and $D=4$, respectively.}
	\label{fit}
\end{figure}

Figure~\ref{LB} illustrates the objective function values over iterations by the adaptive MPMC-IL algorithm. Figure~\ref{fit} shows the approximate marginal posterior distributions derived from the adaptive MPMC-IL algorithm when $D=1$ and $D=4$.  From Figure~\ref{LB} we can see that the objective function value converges after a few iterations and stays an almost flat trend in the following updates. This indicates that the posterior of this model can be fitted well by a Gaussian distribution. Adding more components does not significantly improve the optimization of KL. Despite this, the successive component-updating procedure retains the goodness of fitting, as confirmed by Figure~\ref{fit}. This suggests that the adaptive MPMC-IL algorithm is rather stable and robust for handling uni-modal problems.

As in the g-and-k model, the objective function value in Figure~\ref{LB} first decreases slightly after updating the number of components due to the effect of initialization, and then increases quickly thanks to the main MPMC iterations. Due to the good fitness achieved in the first round of iterations, we do not report the results for the adaptive window size rule here.

\section {Conclusion} \label{sec:con}

In this work, the MPMC algorithm is generalized to the likelihood-free setting.
The proposed algorithm (MPMC-IL) is able to fit complex posteriors with multimodal and heavy tails structures via a proper mixture distribution. Moreover, an adaptive and automatic component-updating algorithm is also provided to drive a better approximation over an extended mixture family. 
The only requirements for the initialial distribution in the adaptive procedure is a Gaussian distribution and the number of components are determined automatically.

For future work, incorporating variance reduction methods in the IS sampling would be a promising avenue, such as the randomized quasi-Monte Carlo methods. We expect that improved sampling accuracy could transfer directly to improved
optimization. In addition, applying the proposed approach to other types of mixture families (such as mixtures of multivariate Student t distributions) would be another interesting topic.

\section*{Acknowledgments}
This work was supported by the National Science Foundation of China (No. 12071154), Guangdong Basic and Applied Basic Research Foundation (No. 2021A1515010275), and Guangzhou Science and Technology Program (No. 202102020407).


\begin{thebibliography}{26}
	\ifx \bisbn   \undefined \def \bisbn  #1{ISBN #1}\fi
	\ifx \binits  \undefined \def \binits#1{#1}\fi
	\ifx \bauthor  \undefined \def \bauthor#1{#1}\fi
	\ifx \batitle  \undefined \def \batitle#1{#1}\fi
	\ifx \bjtitle  \undefined \def \bjtitle#1{#1}\fi
	\ifx \bvolume  \undefined \def \bvolume#1{\textbf{#1}}\fi
	\ifx \byear  \undefined \def \byear#1{#1}\fi
	\ifx \bissue  \undefined \def \bissue#1{#1}\fi
	\ifx \bfpage  \undefined \def \bfpage#1{#1}\fi
	\ifx \blpage  \undefined \def \blpage #1{#1}\fi
	\ifx \burl  \undefined \def \burl#1{\textsf{#1}}\fi
	\ifx \doiurl  \undefined \def \doiurl#1{\url{https://doi.org/#1}}\fi
	\ifx \betal  \undefined \def \betal{\textit{et al.}}\fi
	\ifx \binstitute  \undefined \def \binstitute#1{#1}\fi
	\ifx \binstitutionaled  \undefined \def \binstitutionaled#1{#1}\fi
	\ifx \bctitle  \undefined \def \bctitle#1{#1}\fi
	\ifx \beditor  \undefined \def \beditor#1{#1}\fi
	\ifx \bpublisher  \undefined \def \bpublisher#1{#1}\fi
	\ifx \bbtitle  \undefined \def \bbtitle#1{#1}\fi
	\ifx \bedition  \undefined \def \bedition#1{#1}\fi
	\ifx \bseriesno  \undefined \def \bseriesno#1{#1}\fi
	\ifx \blocation  \undefined \def \blocation#1{#1}\fi
	\ifx \bsertitle  \undefined \def \bsertitle#1{#1}\fi
	\ifx \bsnm \undefined \def \bsnm#1{#1}\fi
	\ifx \bsuffix \undefined \def \bsuffix#1{#1}\fi
	\ifx \bparticle \undefined \def \bparticle#1{#1}\fi
	\ifx \barticle \undefined \def \barticle#1{#1}\fi
	\bibcommenthead
	\ifx \bconfdate \undefined \def \bconfdate #1{#1}\fi
	\ifx \botherref \undefined \def \botherref #1{#1}\fi
	\ifx \url \undefined \def \url#1{\textsf{#1}}\fi
	\ifx \bchapter \undefined \def \bchapter#1{#1}\fi
	\ifx \bbook \undefined \def \bbook#1{#1}\fi
	\ifx \bcomment \undefined \def \bcomment#1{#1}\fi
	\ifx \oauthor \undefined \def \oauthor#1{#1}\fi
	\ifx \citeauthoryear \undefined \def \citeauthoryear#1{#1}\fi
	\ifx \endbibitem  \undefined \def \endbibitem {}\fi
	\ifx \bconflocation  \undefined \def \bconflocation#1{#1}\fi
	\ifx \arxivurl  \undefined \def \arxivurl#1{\textsf{#1}}\fi
	\csname PreBibitemsHook\endcsname
	
	\bibitem{cappe2008adaptive}
	\begin{barticle}
		\bauthor{\bsnm{Capp\'e}, \binits{O.}},
		\bauthor{\bsnm{Douc}, \binits{R.}},
		\bauthor{\bsnm{Guillin}, \binits{A.}},
		\bauthor{\bsnm{Marin}, \binits{J.-M.}},
		\bauthor{\bsnm{Robert}, \binits{C.P.}}:
		\batitle{Adaptive importance sampling in general mixture classes}.
		\bjtitle{Statistics and Computing}
		\bvolume{18}(\bissue{4}),
		\bfpage{447}--\blpage{459}
		(\byear{2008})
	\end{barticle}
	\endbibitem
	
	\bibitem{Min:2013}
	\begin{botherref}
		\oauthor{\bsnm{Minka}, \binits{T.P.}}:
		Expectation propagation for approximate {B}ayesian inference.
		arXiv preprint arXiv:1301.2294
		(2013)
	\end{botherref}
	\endbibitem
	
	\bibitem{Min:2005}
	\begin{botherref}
		\oauthor{\bsnm{Minka}, \binits{T.P.}}:
		Divergence measures and message passing.
		Technical report,
		Citeseer
		(2005)
	\end{botherref}
	\endbibitem
	
	\bibitem{Beau:2010}
	\begin{barticle}
		\bauthor{\bsnm{Beaumont}, \binits{M.A.}}:
		\batitle{Approximate {B}ayesian computation in evolution and ecology}.
		\bjtitle{Annual review of ecology, evolution, and systematics}
		\bvolume{41},
		\bfpage{379}--\blpage{406}
		(\byear{2010})
	\end{barticle}
	\endbibitem
	
	\bibitem{davi:2012}
	\begin{barticle}
		\bauthor{\bsnm{Davison}, \binits{A.C.}},
		\bauthor{\bsnm{Padoan}, \binits{S.A.}},
		\bauthor{\bsnm{Ribatet}, \binits{M.}}:
		\batitle{Statistical modeling of spatial extremes}.
		\bjtitle{Statistical science}
		\bvolume{27}(\bissue{2}),
		\bfpage{161}--\blpage{186}
		(\byear{2012})
	\end{barticle}
	\endbibitem
	
	\bibitem{durbin:2012}
	\begin{bbook}
		\bauthor{\bsnm{Durbin}, \binits{J.}},
		\bauthor{\bsnm{Koopman}, \binits{S.J.}}:
		\bbtitle{Time Series Analysis by State Space Methods},
		\bedition{2}nd edn.
		\bpublisher{Oxford University Press}, \blocation{Oxford}
		(\byear{2012})
	\end{bbook}
	\endbibitem
	
	\bibitem{tavare1997}
	\begin{barticle}
		\bauthor{\bsnm{Tavar{\'e}}, \binits{S.}},
		\bauthor{\bsnm{Balding}, \binits{D.J.}},
		\bauthor{\bsnm{Griffiths}, \binits{R.C.}},
		\bauthor{\bsnm{Donnelly}, \binits{P.}}:
		\batitle{Inferring coalescence times from dna sequence data}.
		\bjtitle{Genetics}
		\bvolume{145}(\bissue{2}),
		\bfpage{505}--\blpage{518}
		(\byear{1997})
	\end{barticle}
	\endbibitem
	
	\bibitem{mar:2003}
	\begin{barticle}
		\bauthor{\bsnm{Marjoram}, \binits{P.}},
		\bauthor{\bsnm{Molitor}, \binits{J.}},
		\bauthor{\bsnm{Plagnol}, \binits{V.}},
		\bauthor{\bsnm{Tavar{\'e}}, \binits{S.}}:
		\batitle{Markov chain {M}onte {C}arlo without likelihoods}.
		\bjtitle{Proceedings of the National Academy of Sciences}
		\bvolume{100}(\bissue{26}),
		\bfpage{15324}--\blpage{15328}
		(\byear{2003})
	\end{barticle}
	\endbibitem
	
	\bibitem{sisson:2018}
	\begin{bbook}
		\bauthor{\bsnm{Sisson}, \binits{S.A.}},
		\bauthor{\bsnm{Fan}, \binits{Y.}},
		\bauthor{\bsnm{Mark}, \binits{B.}}:
		\bbtitle{Handbook of Approximate Bayesian Computation}.
		\bpublisher{Chapman and Hall/CRC}, \blocation{Boca Raton, Florida}
		(\byear{2018})
	\end{bbook}
	\endbibitem
	
	\bibitem{BC:2014}
	\begin{barticle}
		\bauthor{\bsnm{Barthelm{\'e}}, \binits{S.}},
		\bauthor{\bsnm{Chopin}, \binits{N.}}:
		\batitle{Expectation propagation for likelihood-free inference}.
		\bjtitle{Journal of the American Statistical Association}
		\bvolume{109}(\bissue{505}),
		\bfpage{315}--\blpage{333}
		(\byear{2014})
	\end{barticle}
	\endbibitem
	
	\bibitem{tran:2017}
	\begin{barticle}
		\bauthor{\bsnm{Tran}, \binits{M.-N.}},
		\bauthor{\bsnm{Nott}, \binits{D.J.}},
		\bauthor{\bsnm{Kohn}, \binits{R.}}:
		\batitle{Variational {B}ayes with intractable likelihood}.
		\bjtitle{Journal of Computational and Graphical Statistics}
		\bvolume{26}(\bissue{4}),
		\bfpage{873}--\blpage{882}
		(\byear{2017})
	\end{barticle}
	\endbibitem
	
	\bibitem{he:2021}
	\begin{botherref}
		\oauthor{\bsnm{He}, \binits{Z.}},
		\oauthor{\bsnm{Xu}, \binits{Z.}},
		\oauthor{\bsnm{Wang}, \binits{X.}}:
		Unbiased {MLMC}-based variational {B}ayes for likelihood-free inference.
		arXiv preprint arXiv:2109.12728
		(2021)
	\end{botherref}
	\endbibitem
	
	\bibitem{Giles2015}
	\begin{barticle}
		\bauthor{\bsnm{Giles}, \binits{M.B.}}:
		\batitle{Multilevel {M}onte {C}arlo methods}.
		\bjtitle{Acta Numerica}
		\bvolume{24},
		\bfpage{259}--\blpage{328}
		(\byear{2015})
	\end{barticle}
	\endbibitem
	
	\bibitem{Rhee2015}
	\begin{barticle}
		\bauthor{\bsnm{Rhee}, \binits{C.-H.}},
		\bauthor{\bsnm{Glynn}, \binits{P.W.}}:
		\batitle{Unbiased estimation with square root convergence for sde models}.
		\bjtitle{Operations Research}
		\bvolume{63}(\bissue{5}),
		\bfpage{1026}--\blpage{1043}
		(\byear{2015})
	\end{barticle}
	\endbibitem
	
	\bibitem{ong:2018}
	\begin{barticle}
		\bauthor{\bsnm{Ong}, \binits{V.M.}},
		\bauthor{\bsnm{Nott}, \binits{D.J.}},
		\bauthor{\bsnm{Tran}, \binits{M.-N.}},
		\bauthor{\bsnm{Sisson}, \binits{S.A.}},
		\bauthor{\bsnm{Drovandi}, \binits{C.C.}}:
		\batitle{Variational {B}ayes with synthetic likelihood}.
		\bjtitle{Statistics and Computing}
		\bvolume{28}(\bissue{4}),
		\bfpage{971}--\blpage{988}
		(\byear{2018})
	\end{barticle}
	\endbibitem
	
	\bibitem{PBJ:2012}
	\begin{bchapter}
		\bauthor{\bsnm{Paisley}, \binits{J.}},
		\bauthor{\bsnm{Blei}, \binits{D.}},
		\bauthor{\bsnm{Jordan}, \binits{M.}}:
		\bctitle{Variational {B}ayesian inference with stochastic search}.
		In: \bbtitle{Proceedings of the 29th International Coference on International
			Conference on Machine Learning},
		pp. \bfpage{1363}--\blpage{1370}
		(\byear{2012})
	\end{bchapter}
	\endbibitem
	
	\bibitem{MF:2017}
	\begin{bchapter}
		\bauthor{\bsnm{Miller}, \binits{A.C.}},
		\bauthor{\bsnm{Foti}, \binits{N.J.}},
		\bauthor{\bsnm{D'Amour}, \binits{A.}},
		\bauthor{\bsnm{Adams}, \binits{R.P.}}:
		\bctitle{Reducing reparameterization gradient variance}.
		In: \bbtitle{Advances in Neural Information Processing Systems}
		(\byear{2017})
	\end{bchapter}
	\endbibitem
	
	\bibitem{blum_comparative_2013}
	\begin{barticle}
		\bauthor{\bsnm{Blum}, \binits{M.G.B.}},
		\bauthor{\bsnm{Nunes}, \binits{M.A.}},
		\bauthor{\bsnm{Prangle}, \binits{D.}},
		\bauthor{\bsnm{Sisson}, \binits{S.A.}}:
		\batitle{A comparative review of dimension reduction methods in approximate
			{B}ayesian computation}.
		\bjtitle{Statistical Science}
		\bvolume{28}(\bissue{2}),
		\bfpage{189}--\blpage{208}
		(\byear{2013})
	\end{barticle}
	\endbibitem
	
	\bibitem{em:1977}
	\begin{barticle}
		\bauthor{\bsnm{Dempster}, \binits{A.P.}},
		\bauthor{\bsnm{Laird}, \binits{N.M.}},
		\bauthor{\bsnm{Rubin}, \binits{D.B.}}:
		\batitle{Maximum likelihood from incomplete data via the {EM} algorithm}.
		\bjtitle{Journal of the Royal Statistical Society: Series B (Methodological)}
		\bvolume{39}(\bissue{1}),
		\bfpage{1}--\blpage{22}
		(\byear{1977})
	\end{barticle}
	\endbibitem
	
	\bibitem{douc:2007a}
	\begin{barticle}
		\bauthor{\bsnm{Douc}, \binits{R.}},
		\bauthor{\bsnm{Guillin}, \binits{A.}},
		\bauthor{\bsnm{Marin}, \binits{J.-M.}},
		\bauthor{\bsnm{Robert}, \binits{C.P.}}:
		\batitle{Convergence of adaptive mixtures of importance sampling schemes}.
		\bjtitle{The Annals of Statistics}
		\bvolume{35}(\bissue{1}),
		\bfpage{420}--\blpage{448}
		(\byear{2007})
	\end{barticle}
	\endbibitem
	
	\bibitem{douc:2007b}
	\begin{barticle}
		\bauthor{\bsnm{Douc}, \binits{R.}},
		\bauthor{\bsnm{Guillin}, \binits{A.}},
		\bauthor{\bsnm{Marin}, \binits{J.-M.}},
		\bauthor{\bsnm{Robert}, \binits{C.P.}}:
		\batitle{Minimum variance importance sampling via population {M}onte {C}arlo}.
		\bjtitle{ESAIM: Probability and Statistics}
		\bvolume{11},
		\bfpage{427}--\blpage{447}
		(\byear{2007})
	\end{barticle}
	\endbibitem
	
	\bibitem{gunawan2021flexible}
	\begin{botherref}
		\oauthor{\bsnm{Gunawan}, \binits{D.}},
		\oauthor{\bsnm{Kohn}, \binits{R.}},
		\oauthor{\bsnm{Nott}, \binits{D.}}:
		Flexible variational {B}ayes based on a copula of a mixture of normals.
		arXiv preprint arXiv:2106.14392
		(2021)
	\end{botherref}
	\endbibitem
	
	\bibitem{RM:2002}
	\begin{barticle}
		\bauthor{\bsnm{Rayner}, \binits{G.}},
		\bauthor{\bsnm{MacGillivray}, \binits{H.}}:
		\batitle{Weighted quantile-based estimation for a class of transformation
			distributions}.
		\bjtitle{Computational Statistics \& Data Analysis}
		\bvolume{39}(\bissue{4}),
		\bfpage{401}--\blpage{433}
		(\byear{2002})
	\end{barticle}
	\endbibitem
	
	\bibitem{AKM:2009}
	\begin{barticle}
		\bauthor{\bsnm{Allingham}, \binits{D.}},
		\bauthor{\bsnm{King}, \binits{R.A.}},
		\bauthor{\bsnm{Mengersen}, \binits{K.L.}}:
		\batitle{Bayesian estimation of quantile distributions}.
		\bjtitle{Statistics and Computing}
		\bvolume{19}(\bissue{2}),
		\bfpage{189}--\blpage{201}
		(\byear{2009})
	\end{barticle}
	\endbibitem
	
	\bibitem{drovandi:2011}
	\begin{barticle}
		\bauthor{\bsnm{Drovandi}, \binits{C.C.}},
		\bauthor{\bsnm{Pettitt}, \binits{A.N.}},
		\bauthor{\bsnm{Faddy}, \binits{M.J.}}:
		\batitle{Approximate {B}ayesian computation using indirect inference}.
		\bjtitle{Journal of the Royal Statistical Society: Series C (Applied
			Statistics)}
		\bvolume{60}(\bissue{3}),
		\bfpage{317}--\blpage{337}
		(\byear{2011})
	\end{barticle}
	\endbibitem
	
	\bibitem{Fitz:1993}
	\begin{barticle}
		\bauthor{\bsnm{Fitzmaurice}, \binits{G.M.}},
		\bauthor{\bsnm{Laird}, \binits{N.M.}}:
		\batitle{A likelihood-based method for analysing longitudinal binary
			responses}.
		\bjtitle{Biometrika}
		\bvolume{80}(\bissue{1}),
		\bfpage{141}--\blpage{151}
		(\byear{1993})
	\end{barticle}
	\endbibitem
\end{thebibliography}
%

\end{document}